\documentclass[12pt ]{article}

\usepackage{graphicx} 
\usepackage[latin1]{inputenc}
\usepackage{amsmath,amssymb}
\usepackage{latexsym}

\newtheorem{theorem}{Theorem}[section]
\newtheorem{corollary}[theorem]{Corollary}
\newtheorem{lemma}[theorem]{Lemma}

\numberwithin{equation}{section}

\def \be{\begin{equation}}
\def \ee{\end{equation}}
\def \bt{\begin{theorem}}
\def \et{\end{theorem}}
\def \bea{\begin{eqnarray}}
\def \eea{\end{eqnarray}}
\def \bas{\begin{eqnarray*}}
\def \eas{\end{eqnarray*}}

\def \al{\alpha}
\def \bb{\beta}
\def \ga{\gamma}
\def \Ga{\Gamma}
\def \de{\delta}

\def \ep{\epsilon}

\def \la{\lambda}
\def \La{\Lambda}

\def \si{\sigma}

\def \th{\theta}

\def \ff{\infty}
\def \wh{\widehat}
\def \wt{\widetilde}

\def \rar{\rightarrow}

\def \R{{\bf R}}

\def \Z{{\bf Z}}

\def \({\left(}
\def \){\right)}

\def \lc{\left\{}
\def \rc{\right\}}

\def \nn{\nonumber}
\def \Proof{{\bf Proof: }}

\def \bc{\begin{center} }
\def \ec{\end{center} }
\def \bs{\begin{slide} }
\def \es{\end{slide} }

\def\R{{\mathbb R}}
\def\Z{{\mathbb Z}}

\def\phi{\varphi}

\def\proof{{\medskip\noindent {\bf Proof. }}}

\def\qed{{\hfill $\square$ \bigskip}}

\def\Var{{\mathop {{\rm Var\, }}}}

\def\square{{\vcenter{\vbox{\hrule height.3pt
          \hbox{\vrule width.3pt height5pt \kern5pt
             \vrule width.3pt}
          \hrule height.3pt}}}}

 \def\bE {{\Bbb E}}

\def\bP {{\Bbb P}}

\def\square{{\vcenter{\vbox{\hrule height.3pt
           \hbox{\vrule width.3pt height5pt \kern5pt
              \vrule width.3pt}
           \hrule height.3pt}}}}
\def\qed{{\hfill $\square$ \bigskip}}

\def\E{{\bE}}
\def\P{{\bP}}

\begin{document}

\title{Large deviations and renormalization for Riesz potentials of stable intersection measures}
\author{Xia Chen\thanks {Research   partially supported by NSF grant
\#DMS-0704024.}\, \, Jay Rosen\thanks {Research partially supported  by grants
from the NSF  and from PSC-CUNY.}}

\maketitle

\bibliographystyle{amsplain}

\begin{abstract} 
We study  the object formally defined as
\be\gamma\big([0,t]^{2}\big)=\int\!\!\int_{[0,t]^{2}}
\vert X_{s}- X_{r}\vert^{-\sigma}\,dr\,ds-E\int\!\!\int_{[0,t]^{2}}
\vert X_{s}- X_{r}\vert^{-\sigma}\,dr\,ds,\label{jr.1}
\ee
where $X_{t}$ is the  symmetric stable processes of
index $0<\beta\le 2$ in $R^{d}$. When $\beta\le\sigma<\displaystyle\min \Big\{{3\over 2}\beta, d\Big\}$, this has to be defined as a limit, in the spirit of renormalized  self-intersection local time.  We obtain results about the 
large deviations and laws of the iterated logarithm
for $\gamma$. This is applied to obtain results about stable processes in random potentials.
\end{abstract}



\section{Introduction }

Let $X_t$ be a $d$-dimensional 
symmetric stable process of
index $0<\bb\le 2$. Thus we assume that there is a continuous function
$\psi(\lambda) $ on $\R^d $ which is strictly positive for $|\la|\neq 0$, with
$$
\psi(r\lambda)=r^\beta\psi(\lambda)\hskip.1in\hbox{and}\hskip.1in
\psi(-\lambda)=\psi(\lambda),\hskip.2in r>0,\hskip.1in \lambda\in \R^d
$$
such that
\begin{equation}
\E e^{ i\lambda \cdot X_{ t}}=e^{ -t\psi(\lambda)},\hskip.2in t\ge 0,
\hskip.1in
\lambda\in\R^d.\label{a1.0}
\end{equation} 
It follows that there is a constant 
$C >0$
such that
$$
C^{-1}\vert\lambda\vert^\beta\le\psi(\lambda)\le C^{-1}\vert\lambda\vert^\beta,
\hskip.2in \lambda\in\R^d.
$$

In this paper we study
\begin{equation}\label{1.2}
\eta(A)=\int\!\!\int_A
\vert X_s- X_r\vert^{-\sigma}\,dr\,ds,\hskip.2in A\subset R_{+}^{2}
\end{equation} 
and, more generally, 
\begin{equation}
\eta^{z}(A)=\int\!\!\int_{A}
\vert X_s- X_r-z\vert^{-\sigma}\,dr\,ds \label{a1.1z}
\end{equation}
for $z\in \R^{d}$.
We are particularly interested in the case $A=[0, t]^2$ or $[0,t]_<^2$, where for any $t>0$,
$$
[0, t]_<^2=\Big\{(r, s)\in [0, t]^2;\hskip.1in r<s\Big\}.
$$
Thus we will study
\begin{equation}
\eta^{z}([0, t]_<^2)=\int\!\!\int_{[0, t]_<^2}
\vert X_s- X_r-z\vert^{-\sigma}\,dr\,ds.\label{i1.1}
\end{equation}
We can write
\begin{equation}
\eta^{z}([0, t]_<^2)
=\int_{\R^d} {1 \over |x-z|^{\si}}\,\,\mu_{[0, t]_<^2}(\,dx),\label{a1.88}
\end{equation}
where $\mu_{A}$  for  $A\subseteq \R_{+}^{2}$ is the measure on $\R^{d}$  defined by
\begin{equation}
\mu_{A}(B)=\int\int_{A}1_{\{X_s- X_r\in B\}}\,dr\,ds. \label{a1.89}
\end{equation}
We  refer to $\mu_{A}$ as the intersection measure for the stable process $X_{t}$, since whenever $\mu_{[0, t]_<^2}$ has a density $\al_{t}(x)$ which is continuous at $x$, $\al_{t}(x)$ is the intersection local time for $X_{t}$. In particular, if $\al_{t}(x)$  is continuous at 
$x=0$, $\al_{t}(0)$ is a `measure' of the set 
$\{(r, s)\in [0, t]_{<}^2\,|\, X_s=X_r\}$.

(\ref{a1.88}) shows that  $\eta^{z}([0, t]_<^2)$ 
is the Riesz potential of the intersection measure 
$\mu_{[0, t]_<^2}$. 
(In the terminology of \cite{D}, 
$\eta^{z}([0, t]_<^2)$ is 
the Riesz-Frostman potential of the intersection measure.)

Notice that
$\eta\big([0, t]^2)=2\eta\big([0, t]_<^2\big)$.

When $0<\sigma<\min\{\beta, d\}$,
\begin{equation}
\E\int\!\!\int_{[0, t]_<^2}\vert X_s-X_r\vert^{-\sigma}\,dr\,ds
=\E(\vert X_1\vert^{-\sigma})\,\,\int_{0}^{t}\!\!\int_{0}^{s}
{1\over  (s-r)^{\sigma/\beta}}\,dr\,ds\label{p1.1}
\end{equation}
is finite for all $t\geq 0$, so that $\eta\big([0, t]_<^2\big)<\infty$, a.s.

We are interested in Riesz potentials of  intersection measures  for two reasons.
First, our investigation is
motivated by
applications to polymer models. Mathematically,  a
random polymer is modeled as a random path $\omega$ whose probability measure is given in terms of  the Gibbs measure
\begin{align}\label{polymer-1}
P_t(\omega)={1\over Z_t}e^{\pm H_t(\omega)}\,d\omega.
\end{align}
Here $H_t(\omega)\ge 0$ is a suitable  Hamiltonian 
which describes  the  interaction between the monomers along the path
$\omega=\{X_s;\hskip.1in 0\le s\le t\}$, $d\omega$ represents the underlying measure on path space, and $Z_{t}=E\(e^{\pm H_t(\omega)}\)$ is the normalization. In most  models,
the role of $H_t$ is to
reward or penalize   attraction between monomers. The first case, with $+H_t$, describes a ``self-attracting'' polymer, while the second case, with $-H_t$, describes a
 ``self-repelling'' polymer.
We refer to the  recent book by den Hollander, \cite{Hollander}, for a systematic 
overview of  polymer models.

 In the existing literature,
$H_t$ is often taken to be  the self-intersection 
local time, formally defined as 
\begin{align}\label{polymer-2}
H_t=\int_0^t\!\!\int_0^t\delta_0(X_r-X_s)drds.
\end{align}
In this model, the monomers along the path interact only when they intersect.

If one believes that all monomers along the path interact, but the strength of the interaction decreases with distance,  then
the choice of
\begin{align}\label{polymer-3}
H_t=\int_0^t\!\!\int_0^t\vert X_r-X_s\vert^{-\sigma}drds
\end{align}
would be a more realistic model.

The second reason for our interest in Riesz potentials   of  intersection measures     arises from the "polaron problem",
which originated in electrostatics. See \cite{Feynman, Pekar}
 for  general information. The
integral in (\ref{polymer-3}) is associated with
the asymptotics of
the mean-field, or long range interaction, polaron, while the integral
$$
\int_0^t\!\!\int_0^t{e^{-\vert r-s\vert}\over
\vert X_r-X_s\vert^{\sigma}}drds
$$
is associated with a polaron with interactions which are exponentially damped  
in time.
Donsker and Varadhan 
\cite{DV} solved a long standing  problem 
 in physics by showing that, for Brownian motion $W_t$ in $R^{3}$,
\begin{align}\label{D-V}
{\cal D}(\theta)\equiv\lim_{t\to\infty}{1\over t}\log 
\E\bigg\{
\theta\int_0^t\!\!\int_0^t{e^{-\vert r-s\vert}\over\vert W_r-W_s\vert}drds
\bigg\}
\end{align}
exists and
$$
\lim_{\theta\to\infty}{{\cal D}(\theta)\over\theta^2}=
\sup_{g\in{\cal F}_2}\bigg\{\int\!\!\int_{\R^3\times\R^3}
{g^2(x)g^2(y)\over\vert x-y\vert}
dxdy-{1\over 2}\int_{\R^3}\vert\nabla g(x)\vert^2dx\bigg\}
$$
where (with $d=3$ in Donsker and Varadhan's setting)
$$
{\cal F}_2=\big\{g\in {\cal L}^2(\R^d);\hskip.1in \|g\|_2=1\hskip.1in\hbox{and}
\hskip.1in\|\nabla g\|_2<\infty\big\}.
$$

Mansmann (\cite{Mansmann}) showed that for Brownian motion $W_t$ in $R^{3}$, $d\ge 3$,  
\begin{eqnarray} &&
\lim_{t\to\infty}{1\over t}\log \E\exp\bigg\{{1 \over t}  
\int_0^{ t}\!\!\int_0^{ t} {1\over |W_s-W_r|
}\,dr\,ds \bigg\}\label{dv}\\ &&
= \sup_{g\in \mathcal{F}_2}
\bigg\{\int_{\R^d}\int_{\R^d}
{ g^2(x) g^2 (y)\over |x-y|}dx\,dy -
{ 1\over 2}\| \nabla f \|_{ 2}^{ 2}\bigg\}.\nonumber
\end{eqnarray}
He was unable to extend this result to  dimensions $d=1,2$.  

Our first main theorem is the large deviation principle for
$\eta\big([0, t]_<^2\big)$.
For $0<\si<d$ let
\begin{equation}
\phi_{ d-\si}(\lambda)={C_{ d,\si} \over |\la|^{ d-\si}}\label{a1.4k}
\end{equation}
where $C_{ d,\si}=\pi^{-d/2}2^{ -\si}\Ga( { d-\si\over 2})/\Ga( { \si\over
2})$.
 Write
\begin{equation}
\rho=\sup_{\vert\vert f\vert\vert_2=1}
\int_{\R^d}\bigg[\int_{\R^d}{f(\lambda +\gamma)f(\gamma)\over\sqrt{
1+\psi(\lambda +\gamma)}\sqrt{ 1+\psi(\gamma)}}d\gamma\bigg]^2\phi_{
d-\si}(\lambda)d\lambda.\label{a1.7}
\end{equation}
Clearly, $\rho >0$.  According to Lemma 1.6 in \cite{BCR3}, $0<\rho<\infty$ when 
$0<\sigma<\min\{2\beta, d\}$.

 By the scaling property 
 \begin{equation}\label{1.3}
\eta\big([0, t]_<^2\big)\buildrel
d\over =t^{2-\beta^{-1}\sigma}\eta\big([0, 1]_<^2\big), \hskip.2in t\ge 0,
\end{equation}
 we need only consider $\eta\big([0, 1]_<^2\big)$
 in the following theorem.

\bt\label{theo-ldp-1} When $0<\sigma<\min\{\beta, d\}$,
\begin{align}\label{2.2}
\lim_{a\to\infty}a^{-\beta/\sigma}\log\P\Big\{
\eta\big([0,1]_<^2\big)\ge a\Big\}
=-2^{-{\beta/\sigma}}{\sigma\over \beta}
\Big({2\beta -\sigma\over \beta}\Big)^{2\beta -\sigma\over \sigma}
\rho^{-\beta/\sigma}.
\end{align}
\et 
Using scaling, (\ref{1.3}), and  Varadhan's intergral lemma, \cite[Section 4.3]{DZ}, we obtain
the asymptotics 
\begin{align}\label{polymer-4}
&\lim_{t\to\infty}t^{-{2\beta -\sigma\over\beta-\sigma}}\log
\E\exp\bigg\{\int_0^t\!\!\int_0^t\vert X_r-X_s\vert^{-\sigma}drds\bigg\}\\
&\hspace{2 in}={\beta-\sigma\over\beta}\Big({\beta\over 2\beta-\sigma}
\Big)^{2\beta-\sigma\over\beta-\sigma}\rho^{\beta\over\beta-\sigma}\nn
\end{align}
for the partition function in the self-attracting polymer model with
  Hamiltonian defined in (\ref{polymer-3}). 
  Note that when $d=3$, $\beta=2$
and $\sigma=1$ this becomes
 \begin{equation}
 \lim_{t\to\infty}t^{-3}\log
\E\exp\bigg\{\int_0^t\!\!\int_0^t {1 \over \vert W_r-W_s\vert }\,dr\,ds\bigg\}={4 \over 27}\rho^2.\label{ourDV}
 \end{equation} 
Comparison with (\ref{D-V}) shows a striking difference
between the asymptotics of polarons with   long range interactions and those with exponentially damped  interactions.

Theorem \ref{theo-ldp-1}   implies the 
following laws of the iterated
logarithm for 
$\eta\big([0, t]_<^2\big)$.
 \bt\label{theo-lil-1} When $0<\sigma<\min\{\beta, d\}$,
\begin{equation}
\limsup_{t\to\infty}t^{-{2\beta-\sigma\over\beta}}(\log\log t)^{-\sigma/\beta}
\eta\big([0, t]_<^2\big)=2\rho\Big({\beta\over\sigma}\Big)^{\sigma/\beta}
\Big({\beta\over 2\beta -\sigma}\Big)^{ 2\beta -\sigma\over\beta},
\label{2.8}
\end{equation}
almost surely.
\et

We are also  interested in the situation where 
$\beta\le\sigma<\displaystyle\min \Big\{{3\over 2}\beta, d\Big\}$.
In this case $\E\eta\big([0, t]_<^2\big)=\infty$ by (\ref{p1.1}).
We intend to show how to make sense of the object
formally given by
\begin{equation}
\int\!\!\int_{\{0\le r<s\le t\}}\vert X_s-X_r\vert^{-\sigma}drds
-\E\int\!\!\int_{\{0\le r<s\le t\}}\vert X_s-X_r\vert^{-\sigma}drds.\label{1.8}
\end{equation} 

This is reminiscent of the situation for Brownian intersection local time in $R^{2}$. In that case the measure $\mu_{[0, t]_<^2}$ defined in (\ref{a1.89}) has a density $\al_{t}(x)$ which is continuous for all $x\neq 0$, but not for $x=0$. To make sense of $\al_{t}(0)$ we must
 `renormalize'.
This was first done by Varadhan
\cite{rosen:Varadhan}, and has been the subject of a large literature,
see Dynkin \cite{rosen:Dynkin88}, Le Gall \cite{L2},
Bass and Khoshnevisan \cite{rosen:Bass-Khoshnevisan3}, and Rosen
\cite{rosen:RosenJC}. The resulting renormalized
 intersection local time turns out to be the right tool for the
solution of certain ``classical'' problems such as the asymptotic
expansion of the area of the Wiener and stable sausage in the plane
and fluctuations of the range of stable random walks. (See Le Gall
\cite{rosen:LeGall4,rosen:LeGall6}, Le Gall--Rosen
\cite{rosen:LeGall-Rosen} and Rosen \cite{rosen:Rosen92a}).

There are now several ways to `renormalize' 
the Brownian intersection local time $\al_{t}(0)$.       
We briefly recall one such method, since we will use a similar method to renormalize $\eta\big([0, t]_<^2\big)$. Let us write $\al(x,A)$ for the density of the measure $\mu_{A}$. Thus,  $\al_{t}(x)=\al(x,[0, t]_<^2 )$. It can be shown that for any $a<b\leq c<d$, $\mu_{[a,b]\times [c,d]}$ has a continuous density $\al(x,[a,b]\times [c,d])$. We then note that $[0, t]_<^2$ has a decomposition as
\begin{equation}
[0, t]_<^2=\cup_{k=0}^{\ff}\cup_{l=0}^{2^k-1}A_l^k \label{}
\end{equation}
where
\begin{equation}
A_l^k=\Big[{2l\over 2^{k+1}}t, {2l+1\over 2^{k+1}}t\Big)\times
\Big[{2l+1\over 2^{k+1}}t, {2l+2\over 2^{k+1}}t\Big).\label{1.9}
\end{equation} 
It can then be shown that
\begin{equation}
\sum_{k=0}^{\ff}\sum_{l=0}^{2^k-1}\(\al(0,A_l^k)-E(\al(0,A_l^k))\)\label{i1.9}
\end{equation} 
converges, and the limit  is called renormalized
 intersection local time.
 
In the following section we shall carry out a similar 
program to renormalize  $\eta\big([0, t]_<^2\big)$, which 
will give meaning to the formal expression in (\ref{1.8}). 
The resulting object will be denoted by $\gamma\big([0,t]_<^2\big)$. 
This will  make perfectly good sense in the context of  (\ref{polymer-1}), which will now 
  be ``renormalized'' 
to
\begin{align}\label{polymer-5}
P_t(\omega)={1\over \widetilde{Z}_t}e^{\pm\gamma([0, t]^2)}\,d\omega.
\end{align}

The second main result of this paper  is to show that 
 $\gamma\big([0,t]_<^2\big)$ has large deviation properties 
and laws of the iterated
logarithm similar to those established above for 
$\eta\big([0, t]_<^2\big)$ when $0<\sigma<\min\{\beta, d\}$.

\begin{theorem}\label{theo-ldp-2} 
When $\displaystyle 0<\sigma<\min\Big\{{3\over 2}\beta, d\Big\}$,
\begin{align}\label{2.3}
\lim_{a\to\infty}a^{-\beta/\sigma}\log\P\Big\{
\gamma\big([0,1]_<^2\big)\ge a\Big\}
=-2^{-{\beta/\sigma}}{\sigma\over \beta}
\Big({2\beta -\sigma\over \beta}\Big)^{2\beta -\sigma\over \sigma}
\rho^{-\beta/\sigma}.
\end{align}
\end{theorem}

Consider the special case $\beta =\sigma<d$.
Combined with the scaling property given in (\ref{gamma-1})
below, Theorem \ref{theo-ldp-2} shows that the self-attracting polymer measure
(\ref{polymer-5}) `collapses' in   finite time,
by which we mean that 
\begin{align}\label{polymer-6}
\E\exp\Big\{\gamma\big([0, t]^2\big)\Big\}
\left\{\begin{array}{ll}\displaystyle <\infty\hskip.3in t<\rho^{-1}\\
\displaystyle=\infty\hskip.3in t>\rho^{-1}.\end{array}\right.
\end{align}

Theorem \ref{theo-ldp-2} also indicates that the self-attracting polymer
with $H_t=\gamma\big([0, t]_<^2\big)$ is not well defined when $\sigma>\bb$. 
But it is not hard to show that  
$$
\E\exp\Big\{-\gamma\big([0, t]_<^2\Big\}<\infty
$$
for every $t>0$ if
$\displaystyle 0<\sigma<\min\Big\{{3\over 2}\beta, d\Big\}$. A problem relevant to the self-repelling polymer is
to investigate the asymptotics of
$$
\E\exp\Big\{-\eta\big([0, t]_<^2\big)\Big\}\hskip.1in\hbox{or}\hskip.1in
\E\exp\Big\{-\gamma\big([0, t]_<^2\big)\Big\}\hskip.2in 
$$
as $t\to\infty$.
We leave this to future study.

Theorem \ref{theo-ldp-2}  implies the 
following laws of the iterated
logarithm for 
$\ga\big([0, t]_<^2\big)$.
 \bt\label{theo-lil-2} When   $\displaystyle\beta\le \sigma<\min\Big\{{3\over 2}\beta, d\Big\}$,
 \begin{equation}
\limsup_{t\to\infty}t^{-{2\beta-\sigma\over\beta}}(\log\log t)^{-\sigma/\beta}
\gamma\big([0, t]_<^2\big)=2\rho\Big({\beta\over\sigma}\Big)^{\sigma/\beta}
\Big({\beta\over 2\beta -\sigma}\Big)^{ 2\beta -\sigma\over\beta}, \label{2.10}
 \end{equation}
 almost surely.
\et

We next describe an application  to the study of stable processes in a
Brownian potential. 
Let $W$ denote white noise on $L^{2}(R^{d},\,dx)$. That is, 
for every $f\in L^{2}(R^{d},\,dx)$, 
$W(f)$ is the mean zero Gaussian process  with covariance
\begin{equation}
E(W(f)W(g))=\int_{R^{d}}f(x)g(x)\,dx,\label{f1}
\end{equation}
which we take to be independent of our stable process $X_{t}$.
$W(f)$ can be considered as a stochastic integral
\begin{equation}
W(f)=\int_{R^{d}}f(x)\,W(dx)\label{f2}
\end{equation}
with respect to the Brownian sheet, \cite{W}.

Recall the identity
\begin{align}\label{eta-3}
\int_{\R^d}\vert x-z\vert^{-{\sigma +d\over 2}}\vert y-z\vert^{-{\sigma +d\over 2}}
dz=C{1\over\vert x-y\vert^\sigma},\hskip.2in x, y\in\R^d
\end{align}
where 
\begin{align}\label{2.6'}
C=\pi^{d/2}{\displaystyle\Gamma^{2}\Big({d-\sigma\over 4}\Big)
\Gamma\Big({\sigma\over 2}\Big)\over 
\displaystyle\Gamma^{2}\Big({d+\sigma\over 4}\Big)
\Gamma\Big({d-\sigma\over 2}\Big)},
\end{align}
see \cite[p. 118,158]{D} or \cite[p. 118. (8)]{S}.
Then, if we set 
\be
\xi(t, x)=\int_0^t\vert X_s-x\vert^{-{\sigma +d\over 2}}\,ds,\label{eta-3a}
\ee
we see  by Fubini's theorem that
\begin{align}\label{eta-4}
\eta\big([0, t]^2\big)=C\int_{\R^d}\xi(t, x)^2 dx<\ff, \hspace{.2 in}a.s.
\end{align}
 
Hence, almost surely with respect to $X$,
\begin{equation}
F(t)\equiv W(\xi(t,\cdot))=\int_{\R^d}\xi(t,x)W(dx)\label{i1.7}
\end{equation}
is a   mean zero normal random variable with 
\begin{equation}
E\(F^{2}(t)\)={2\over C}\eta\big([0, t]^2_<\big).\label{i1.8}
\end{equation}

By a stable process in a
Brownian potential we mean  the process described by the measure 
\begin{align}\label{spbp-1}
Q_t ={1\over N_t}e^{-\int_0^t V(X_s)ds}\,P_{t},
\end{align}
where $P_{t}$ is the probability for the stable process
$ \{X_s;\hskip.1in 0\le s\le t\}$,   $Z_{t} $ is the normalization
and  the potential function $V(x)$ is of the form
\be
V(x)=\int_{\R^d}K(y-x)W(dy),\label{spbp-2}
\ee
where $K(x)$ is a function on $\R^d$ known as the 
shape function. Roughly speaking,  (\ref{spbp-2})
represents the interaction 
with a field generated by a cloud of
electrons $W(dy)$ with random signed  charges and 
  locations in $\R^d$. 
We refer \cite{CM} and \cite{Sznitman} 
for more  information.

 When $K(x)=\delta_0(x)$, the action $\int_0^t V(X_s)ds$ corresponds to a stable process in Brownian 
scenery (\cite{KS}, \cite{KhoshnevisanLewis}, \cite{Bolthausen}
and \cite{CK}). When $K(x)$ is bounded and locally supported, the long term 
behavior of $\int_0^t V(X_s)ds$  is similar to that of the stable process in Brownian
scenery.

Let $\displaystyle{d\over 2}<p<\min\Big\{d, {d+\beta\over 2}\Big\}$.
The random potential function
\begin{align}\label{potential-2}
V_p(x)\equiv\int_{\R^d}{1\over\vert y-x\vert^p}W(dy),\hskip.2in x\in\R^d
\end{align}
has a very intuitive physical meaning. 
When $d=3$ and $p=1$,   it represents the 
electrostatic potential energy generated by a cloud of
electrons $W(dy)$.
Unfortunately,  the random potential function (\ref{potential-2}) 
is not well-defined,  since $\E V_p^2(x)=\infty$ for every
$x\in\R^d$. However, as we have seen,
\begin{align}\label{i1.7}
F(t)\equiv\int_{\R^d}\bigg[\int_0^t{ds\over\vert y-X_s\vert^p}\bigg]W(dy)
\end{align}
is well defined. Here we have taken $\sigma=2p-d$.  Because of (\ref{spbp-1}), we refer to $F(t)$ as the action for a stable process in a
Brownian potential.

It is easy to see that for each $t>0$,
\begin{align}\label{2.4}
F(t)\buildrel d\over =t^{2\beta-2p+d\over 2\beta}F(1).
\end{align}
The following corollaries  about large deviations and laws of the iterated logarithm for $F(t)$
will follow from Theorem \ref{theo-ldp-1}.
\begin{corollary}\label{theo-ldp-f} If $\displaystyle
{d\over 2}<p<\min\Big\{d, {d+\beta\over 2}\Big\}$, then
\begin{align}\label{2.6}
&\lim_{a\to\infty}a^{-{2\beta\over 2p-d +\beta}}
\log\P\big\{\pm F(1)\ge a\big\}\\
&=-{\beta+2p-d\over\beta}
\Big({C_p\over 8\rho_p}\Big)^{\beta\over \beta+2p-d}
\Big({2\beta-2p+d\over\beta}\Big)^{2\beta-2p-d\over \beta+2p-d}\nn
\end{align}
where the constants  
$$
C_p=\pi^{d/2}{\displaystyle\Gamma^{2}\Big({d-p\over 2}\Big)
\Gamma\Big({2p-d\over 2}\Big)\over 
\displaystyle\Gamma^{2}\Big({p\over 2}\Big)
\Gamma({d-p})},
$$
and
$$
\rho_p=\sup_{\vert\vert f\vert\vert_2=1}
\int_{\R^d}\bigg[\int_{\R^d}{f(\lambda +\gamma)f(\gamma)\over\sqrt{
1+\psi(\lambda +\gamma)}\sqrt{ 1+\psi(\gamma)}}\,d\gamma\bigg]^2
\phi_{2(d-p)}(\lambda)\,d\lambda
$$
come from $C$ and $\rho$ defined in (\ref{2.6'}) and (\ref{a1.7}), 
respectively,
with $\sigma=2p-d$.
\end{corollary}

\begin{corollary}\label{theo-lil-f}  
\begin{align}\label{2.9}
&\limsup_{t\to\infty}t^{-{2\beta-2p-d\over 2\beta}}
(\log\log t)^{-{2p-d+\beta\over 2\beta}}\Big\{\pm F(t)\Big\}\\
&=\sqrt{8\rho_p\over C_p}\Big({\beta\over2p-d+\beta}\Big)^{2p-d+\beta\over 2\beta}
\Big({\beta\over 2\beta-2p+d}\Big)^{2\beta-2p+d\over 2\beta},
\hskip.2in a.s.\nn
\end{align}

\end{corollary}

We have also obtained a variational expression for $\rho.$   Let  
\begin{equation}
 \mathcal{E}_{ \bb}( f,f)=:(2\pi)^{-d}\int_{ \R^d} |\la|^{ \bb}|\wh{f}( \la)|^2\,d\la,
\label{7.0a}
\end{equation} and
\begin{equation}
\mathcal{F}_{ \bb}=\{f\in L^2( R^d)\,|\,\|f\|_{ 2}=1\,,\,
\,\mathcal{E}_{ \bb}( f,f)<\ff\}.\label{7.0b}
\end{equation}
We show in \cite{BCR3} that 
\begin{equation}
\La_{\si}=:\sup_{g\in{\cal F}_\beta}\bigg\{ 
\bigg(\int_{ (\R^{ d})^{ 2}} { g^{ 2}( x )g^{ 2}( y )\over |x+y|^{ \si}}
\,d x\,\,dy\bigg)^{1/2}-{\cal E}_\bb (g,g)\bigg\}<\ff\label{67.0v}
\end{equation} 
when $0<\sigma<\min\{2\beta, d\}$ and we derive a relation, \cite[(1.20)]{BCR3}, between $\rho$ and $ \La_{\si}$.

Outline: In Section \ref{sec-renorm} 
we show how to renormalize $ \eta ([0, t]_<^2 )$ 
when $\beta\le\sigma<\displaystyle\min  (3 \beta/2, d) $, 
and establish 
some exponential estimates for $ \eta ([0, t]_<^2 )$ and  
$ \ga ([0, t]_<^2 )$. In Section \ref{sec-moment asymptotics} 
we establish high moment asymptotics for 
 $ \eta ([0, t]_<^2 )$ and  $ \ga ([0, t]_<^2 )$ which 
lead to the proof of our main theorems on large deviations, 
Theorems \ref{theo-ldp-1} 
and \ref{theo-ldp-2}, in Section \ref{sec-largedev}. 
The corressponding laws of the iterated logarithm are 
established in  Section \ref{sec-lil}. Section \ref{sec-F} 
deals with $F(t)$, the action for a stable process in a random potential. 
Finally, in a short Appendix, Section \ref{sec-Mlim}, 
we give details on approximation of $\rho$ which is 
used in Section \ref{sec-moment asymptotics}.

Conventions: We define
\begin{equation}
\wh{f}( \la)=\int_{\R^d}e^{ ix\cdot\la}f( x)\,dx.\label{c6.1}
\end{equation}
With this notation
\begin{equation}
f( x)=( 2\pi)^{ -d}\int_{\R^d}e^{ -ix\cdot\la}\wh{f}( \la)\,dx,\label{c6.2}
\end{equation}
\begin{equation}
\wh{f\ast g}( \la)=\wh{f}( \la)\wh{g}( \la),\hspace{ .4in}\wh{f g}( \la)
=( 2\pi)^{ -d}\wh{f}( \la)\ast
\wh{g}(\la),\label{c6.3}
\end{equation}
and Parseval's identity is 
\begin{equation}
\langle f,g\rangle=( 2\pi)^{ -d}\langle\wh{f}, \wh{g}\rangle.\label{c6.4}
\end{equation}
If $\Phi\in \mathcal{S'}(\R^{d})$, the set of tempered distributions on $\R^{d}$, we use 
$\mathcal{F}(\Phi)$ to denote the Fourier transform of $\Phi$, so that for any 
$f\in \mathcal{S}(\R^{d})$
\begin{equation}
\mathcal{F}(\Phi)(f)=\Phi (\wh{f}).\label{a1.4j}
\end{equation}
It is well known. e.g. \cite[p. 156]{D}, that 
$\phi_{ d-\si}\in \mathcal{S'}(\R^{d})$ for any $0<\si<d$  
 and 
\begin{equation}
\mathcal{F}(\phi_{ d-\si})={1 \over |x|^{\si}}.\label{a1.4}
\end{equation}

\section{Renormalization }\label{sec-renorm}

We begin by proving an exponential integrability result for $\eta\big([0, 1]_<^2\big)$.

\begin{theorem}\label{eta-1} If $0<\sigma<\min\{\beta, d\}$, there
is a $c>0$ such that
\begin{align}\label{eta-2}
\E\exp\Big\{c\eta\big([0, 1]_<^2\big)^{\beta/\sigma}\Big\}<\infty.
\end{align}
\end{theorem}

\proof 
Recall that by (\ref{eta-4}), for each $t>0$, $\xi(t, x)\in L^{2}(R^{d},\,dx)$   almost surely. 
By the  triangle inequality, for any $s, t>0$,
$$
\begin{aligned}
&\bigg\{\int_{\R^d}\xi(s+t, x)^2 dx\bigg\}^{1/2}\\
&\le \bigg\{\int_{\R^d}\xi(s, x)^2 dx\bigg\}^{1/2}
+\bigg\{\int_{\R^d}\big[\xi(s+t, x)-\xi(s, x)\big]^2 dx\bigg\}^{1/2}.
\end{aligned}
$$
Notice that the integral
\bea
&&
\int_{\R^d}\big[\xi(s+t, x)-\xi(s, x)\big]^2 dx\nonumber\\
&&=C^{-1}\int\!\!\int_{[s,s+t]^{2}}\vert X_u-X_v\vert^{-\sigma}\,du\,dv\nonumber\\
&&=C^{-1}\int\!\!\int_{[0,t]^{2}}\vert X_{s+u}-X_{s+v}\vert^{-\sigma}\,du\,dv
\eea
is independent of $\{X_u;\hskip.1in 0\le u\le s\}$ and has the same
distribution as
$$
\int_{\R^d}\xi(t, x)^2 dx.
$$
The process
\be
\bigg\{\int_{\R^d}\xi(t, x)^2 dx\bigg\}^{1/2},\hskip.2in t\ge 0,\label{p1.5}
\ee
is therefore   sub-additive. By Theorem 1.3.5 in \cite{Chen},
$$
\E\exp\bigg\{\theta\bigg\{\int_{\R^d}\xi(t, x)^2 dx\bigg\}^{1/2}\bigg\}<\infty,
\hskip.2in \forall \theta, t>0,
$$
and for any $\theta>0$, the limit
$$
L(\theta)\equiv \lim_{t\to\infty}{1\over t}\log
\E\exp\bigg\{\theta\bigg\{\int_{\R^d}\xi(t, x)^2 dx\bigg\}^{1/2}\bigg\}
$$
exists with $0\le L(\theta)<\infty$. Using Chebyshev's
inequality we find that
$$
\limsup_{t\to\infty}{1\over t}\log\P\bigg\{
\int_{\R^d}\xi(t, x)^2 dx\ge t^2\bigg\}\le -l
$$
for some $l>0$. By the scaling property given in (\ref{1.3}),
$$
\limsup_{t\to\infty}{1\over t}\log\P\bigg\{
\int_{\R^d}\xi(1, x)^2 dx\ge t^{\sigma/\beta}\bigg\}\le -l,
$$
which leads to (\ref{eta-2}). \qed

We now show how to renormalize $\eta\big([0, t]_<^2\big)$ when 
$\beta\le\sigma<\displaystyle\min \Big\{{3\over 2}\beta, d\Big\}$.
Recall that in this case $\E\eta\big([0, t]_<^2\big)=\infty$.
We will show how to make sense of the object
formally given by (\ref{1.8}). 

 To proceed further, let $\widetilde{X}_t$ be an independent
copy of $X_t$ and define the random measure
\begin{equation}\label{1.4}
\zeta(A)=\int\!\!\int_A\vert X_s- \widetilde{X}_t\vert^{-\sigma}dsdt
\hskip.2in A\subset (\R^+)^2.
\end{equation} 
By  \cite[Theorem 1.1]{BCR3}, $\zeta(A)<\infty$ a.s. for every bounded $A$ and
\begin{equation}\label{1.5}
\zeta \big([0,t]^2\big)\buildrel d \over =
t^{2-\sigma/\beta }\zeta\big([0, 1]^2\big), \hskip.2in t\ge 0.
\end{equation} 
Further, by   \cite[Theorem 1.2]{BCR3}   there is a $\theta>0$ such that
\begin{equation}\label{1.6}
\E\exp\Big\{\theta \zeta\big([0, 1]^2\big)^{\beta/\sigma}\Big\}<\infty.
\end{equation} 
Note that for any $0\le a<b<c<\infty$,
\begin{equation}\label{1.7}
\eta \big([a, b]\times [b, c]\big)\buildrel d\over =
\zeta \big([0,b-a]\times[0, c-b]\big).
\end{equation}

To make sense of (\ref{1.8}) we shall use
 Varadhan's triangular approximation (see, e.g.,
Proposition 6, p.194, \cite{L2}).
Let $t>0$ be fixed.
Recall the subsets $A_l^k$ defined in (\ref{1.9}). By (\ref{1.5}) and (\ref{1.7}) we have  
\begin{equation}\label{1.10}
\eta (A_l^k)\buildrel d\over =\zeta\Big(\Big[0,\hskip.05in
{t\over 2^{k+1}}\Big]^2\Big)
\buildrel d\over =2^{-(k+1)(2-\sigma/\beta)}\zeta\big([0,t]^2).
\end{equation}
In addition, for each $k\ge 0$, the finite sequence
$$
\eta (A_l^k),\hskip.2in l=0,1,\cdots, 2^k-1
$$
is independent. Consequently, 
$$
\Var\bigg(\sum_{l=0}^{2^k-1}\eta (A_l^k)\bigg)
=C2^{-(3-2\sigma/\beta)k}.
$$
Under $0<\sigma<\min\{3\beta/2, d\}$, therefore,
$$
\begin{aligned}
&\bigg\{\E\bigg[\sum_{k=0}^\infty\bigg\vert\sum_{l=0}^{2^k-1}\Big(\eta(A_l^k)
-\E\eta (A_l^k)\Big)\bigg\vert\bigg]^2\bigg\}^{1/2}\\
&\le\sum_{k=0}^\infty\bigg\{\Var\bigg(\sum_{l=0}^{2^k-1}\eta (A_l^k)\bigg)
\bigg\}^{1/2}<\infty.
\end{aligned}
$$
Consequently, the random series
$$
\sum_{k=0}^\infty\bigg\{\sum_{l=0}^{2^k-1}\Big(\eta(A_l^k)
-\E\eta (A_l^k)\Big)\bigg\}
$$
convergences in $L^2(\Omega, {\cal A}, \P)$.
We may therefore  define
$$
\gamma\big([0, t]_<^2\big)
=\sum_{k=0}^\infty\bigg\{\sum_{l=0}^{2^k-1}\Big(\eta(A_l^k)
-\E\eta (A_l^k)\Big)\bigg\}.
$$
This will be our definition of the the object formally  given in (\ref{1.8}).
As in (\ref{1.3}), 
\begin{equation}\label{gamma-1}
\gamma\big([0, t]_<^2\big)\buildrel d\over =t^{2- \sigma/\beta}
\gamma\big([0, 1]_<^2\big),\hskip.2in t\ge 0.
\end{equation}
As in  Theorem \ref{eta-1}, we have the following exponential integrability.

\begin{theorem}\label{1.11} If
$0<\sigma<\displaystyle\min\Big\{{3\over 2}\beta, d\Big\}$, there is a 
$c>0$
such that
\begin{equation}\label{1.12}
\E\exp\Big\{c\big\vert\gamma\big([0, 1]_<^2\big)\big\vert^{\beta/\sigma}
\Big\}<\infty.
\end{equation}
\end{theorem}

\proof 
When $\si<\beta$ this follows  trivially from  Theorem \ref{eta-1}. We can therefore assume that 
$p=\beta/\sigma\leq 1$.
Again, we use the triangular
approximation based on the partition (\ref{1.9})
 with
$t=1$.

In view of (\ref{1.6}) and (\ref{1.10}), applying Lemma 1, \cite{BCR}
to the family of the i.i.d. sequences
$$
\Big\{2^{(k+1)(2- \sigma/\beta)}\big[\zeta(A_l^k)-\E\zeta(A_l^k)\big];
\hskip.1in l=0, 1,\cdots, 2^k-1\Big\},
\hskip.2in k=0, 1,\cdots
$$
gives that for some $\theta>0$
$$
\sup_{k\ge 0}\E\exp\bigg\{\theta \Big\vert 2^{-k/2}\sum_{l=0}^{2^k-1}
2^{k(2- \sigma/\beta)}\big[\zeta(A_l^k)-\E\zeta(A_l^k)\big]
\Big\vert^{\beta/\sigma}\bigg\}<\infty,
$$
or
\begin{align}\label{1.13}
e^{C}\equiv\sup_{k\ge 0}\E\exp\bigg\{2^{ak}\theta \Big\vert\sum_{l=0}^{2^k-1}
\big[\zeta(A_l^k)-\E\zeta(A_l^k)\big]
\Big\vert^{\beta/\sigma}\bigg\}<\infty
\end{align}
where
$$
a={3\beta\over 2\sigma}-1>0.
$$

For each $N\ge 1$ set
$$
b_1=\theta,\hskip.1in b_N=\theta\prod_{j=2}^N\Big(1-2^{-a(j-1)}\Big),
\hskip.2in N=2,3,\cdots.
$$

By H\"older's inequality and the triangular inequality  
$$
\begin{aligned}
&\E\exp\bigg\{
b_N\bigg\vert\sum_{k=0}^N\sum_{l=0}^{2^k-1}\big[\beta (A_l^k)
-\E\beta (A_l^k)\big]\bigg\vert^{\beta/\sigma}\bigg\}\\
&\le\Bigg[\E\exp\bigg\{
b_{N-1}\bigg\vert\sum_{k=0}^{N-1}\sum_{l=0}^{2^k-1}\big[\beta (A_l^k)
-\E\beta (A_l^k)\big]\bigg\vert^{\beta/\sigma}\bigg\}\Bigg]^{1-2^{-a(N-1)}}\\
&\times\Bigg[\E\exp\bigg\{
2^{a(N-1)}b_N\bigg\vert\sum_{l=0}^{2^N-1}\big[\beta (A_l^N)
-\E\beta (A_l^N)\big]\bigg\vert^{\beta/\sigma}\bigg\}\Bigg]^{2^{-a(N-1)}}.
\end{aligned}
$$
Notice that $b_N\le \theta$. By (\ref{1.13}) we have
$$
\begin{aligned}
&\E\exp\bigg\{
b_N\bigg\vert\sum_{k=0}^N\sum_{l=0}^{2^k-1}\big[\beta (A_l^k)
-\E\beta (A_l^k)\big]\bigg\vert^{\beta/\sigma}\bigg\}\\
&\le \exp\Big\{C2^{-a(N-1)}\Big\}\E\exp\bigg\{
b_{N-1}\bigg\vert\sum_{k=0}^{N-1}\sum_{l=0}^{2^k-1}\big[\beta (A_l^k)
-\E\beta (A_l^k)\big]\bigg\vert^{\beta/\sigma}\bigg\}.
\end{aligned}
$$

Repeating the above procedure gives
$$
\begin{aligned}
&\E\exp\bigg\{
b_N\bigg\vert\sum_{k=0}^N\sum_{l=0}^{2^k-1}\big[\beta (A_l^k)
-\E\beta (A_l^k)\big]\bigg\vert^{\beta/\sigma}\bigg\}\\
&\le\exp\Big\{C\sum_{k=0}^N2^{-a(k-1)}\Big\}
\le\exp\Big\{C\Big(1-2^{-a}\Big)^{-1}\Big\}<\infty.
\end{aligned}
$$
Observe that 
$$
b_\infty =
\theta\prod_{j=2}^\infty\Big(1-2^{-a(j-1)}\Big)>0.
$$
By Fatou's lemma, letting $N\to\infty$  we have
$$
\E\exp\Big\{b_\infty\big\vert\gamma \big([0,1]_<^2\big)\big\vert^{\beta/\sigma}
\Big\}
\le\exp\Big\{C'\Big(1-2^{-a}\Big)^{-1}\Big\}<\infty.
$$
\qed

Among other things, we now show that the family
$$
\Big\{\gamma \big([0,t]_<^2\big)\hskip.1in t\ge 0\Big\}
$$
has a continuous version.

\begin{lemma}\label{renorm-10} Assume
$0<\sigma<\displaystyle\min\Big\{{3\over 2}\beta, d\Big\}$.
For any $T>0$ there is a $c=c(T)>0$ such that
$$
\sup_{\scriptstyle s,t\in [0, T]
\atop\scriptstyle
s\not = t}\E\exp\bigg\{c{\big\vert \gamma \big([0,t]_<^2\big)
- \gamma \big([0,s]_<^2\big)\big\vert^{\beta/\sigma}\over
\vert t-s\vert^{(2\beta -\sigma)/(2\sigma)}}\bigg\}
<\infty.
$$
\end{lemma}

\proof For $0\le s<t\le T$,
$$
\gamma \big([0,t]_<^2\big)
- \gamma \big([0,s]_<^2\big)
=\gamma \big([s,t]_<^2\big)+\gamma \big([0,s]\times [s,t]\big).
$$
Notice that
$$
\gamma \big([s,t]_<^2\big)\buildrel d\over =(t-s)^{2- \sigma/\beta}
\gamma([0, 1]_<^2\big).
$$
By Theorem \ref{1.11}, there is a $c_1>0$
such that
$$
\sup_{\scriptstyle s,t\in [0, T]
\atop\scriptstyle
s<t}\E\exp\bigg\{c_1{\big\vert \gamma \big([s,t]_<^2\big)
\big\vert^{\beta/\sigma}\over\vert t-s\vert^{(2\beta-\sigma)/\sigma -1}}\bigg\}
<\infty.
$$
In addition, by  (\ref{1.7})
$$
\gamma \big([0,s]\times [s,t]\big)
\buildrel d\over=\zeta\big([0, s]\times [0, t-s]\big)
-\E\zeta\big([0, s]\times [0, t-s]\big).
$$
In view of (\ref{eta-3}),
$$
\zeta\big([0, s]\times [0, t-s]\big)=C\int_{\R^d}\xi(s, x)\tilde{\xi}(t-s, x)dx,
$$
where
$$
\xi(t,x)=\int_0^t\vert X_u-x\vert^{-{\sigma+d\over 2}}du\hskip.1in\hbox{and}
\hskip.1in\tilde{\xi}(t,x)=\int_0^{t}\vert \widetilde{X}_u-x\vert^{-{\sigma+d\over 2}}du.
$$
By independence, for any integer $m\ge 1$,
$$
\begin{aligned}
&\E\Big[\zeta\big([0, s]\times [0, t-s]\big)^m\Big]\\
&=C^m\int_{(\R^d)^m}dx_1\cdots dx_m\bigg[\E\prod_{k=1}^m\xi(s,x_k)\bigg]
\bigg[\E\prod_{k=1}^m\xi(t-s,x_k)\bigg]\\
&\le C^m\bigg\{\int_{(\R^d)^m}dx_1\cdots dx_m
\bigg[\E\prod_{k=1}^m\xi(s,x_k)\bigg]^2
\bigg\}^{1/2}\\
&\times\bigg\{\int_{(\R^d)^m}dx_1\cdots dx_m\bigg[\E\prod_{k=1}^m\xi(t-s,x_k)\bigg]^2\bigg\}^{1/2}\\
&=\Big\{\E\Big[\zeta\big([0, s]^2\big)^m\Big\}^{1/2}
\Big\{\E\Big[\zeta\big([0, t-s]^2\big)^m\Big\}^{1/2}\\
&=
s^{{2\beta -\sigma\over 2\beta}m}(t-s)^{{2\beta -\sigma\over 2\beta}m}
\E\Big[\zeta\big([0, 1]^2\big)^m\Big\}.
\end{aligned}
$$

Checking (\ref{1.6}), there is $c_2>0$ such that
$$
\sup_{\scriptstyle s,t\in [0, T]
\atop\scriptstyle
s<t}\E\exp\bigg\{c_2{\big\vert \gamma \big([0, s]\times[s, t]\big)
\big\vert^{\beta/\sigma}\over\vert t-s\vert^{2\beta -\sigma\over 2\sigma}}\bigg\}
<\infty.
$$

The proof is now complete. \qed

Recall that a function $\Psi$: $\R^+\longrightarrow\R^+$ is a Young function
if it is convex,
increasing and satisfies $\Psi(0)=0$, $\lim_{x\to\infty}\Psi(x)=\infty$.
The Orlicz space ${\cal L}_\Psi(\Omega, {\cal A}, \P)$
is defined as the linear space of all random variables $X$
on the probability
space $(\Omega, {\cal A}, \P)$ such that
$$
\|X\|_\Psi =\inf\{c>0;\hskip.1in \E\Psi(c^{-1}\vert X\vert)\le 1\}<\infty
$$
It is known that $\|\cdot\|_\Psi$ defines a norm 
(called Orlicz norm) under which ${\cal L}_\Psi(\Omega, {\cal A}, \P)$
becomes a Banach space.

We now choose the Young function $\Psi(\cdot)$ such that
$\Psi(x)\sim \exp\{x^{\beta/\sigma}\}$ as $x\to\infty$.
By Lemma \ref{renorm-10}, for any $T>0$
 there is $c=c(d, \beta, \sigma, T)>0$ such that
 \begin{equation}
 \big\|\gamma \big([0,t]_<^2\big)- \gamma \big([0,s]_<^2\big)\big\|_\Psi
\le c\vert t-s\vert^{(2\beta -\sigma)/(2\beta)},\hskip.2in s, t\in [0, T].\label{jr7}
 \end{equation} 
By Lemma 9 in \cite{CLR}, the family
$\Big\{\gamma \big([0,t]_<^2\big)\hskip.1in t\ge 0\Big\}$ has a continuous
version and in the future we will always use this version. 
Furthermore, for any $0<b<{2\beta -\sigma\over 2\beta}$,
$$
\bigg\|\sup_{0\le s<t\le 1}
{\gamma \big([0,t]_<^2\big)- \gamma \big([0,s]_<^2\big)\over
\vert t-s\vert^b}\bigg\|_\Psi<\infty.
$$
Equivalently, there is
$c>0$ such that
\begin{align}\label{1.14}
\E\exp\bigg\{c\sup_{0\le s<t\le 1}{\big\vert
\gamma\big([0, t]_<^2\big)-\gamma\big([0, s]_<^2\big)\big\vert^{\beta/\sigma}\over
\vert t-s\vert^{b\beta/\sigma}}\bigg\}<\infty.
\end{align}

\section{High moment asymptotics for 
the smoothed version}\label{sec-moment asymptotics}

Throughout this section we assume that 
\begin{equation}
0<\sigma<\min\{2\beta, d\}.\label{p3.10}
\end{equation}
We define the probability density function   
\begin{equation}
h(x)=C^{-1}\prod_{j=1}^d\Big({2\sin x_j\over x_j}\Big)^2,
\hskip.2in x=(x_1,\cdots, x_d)\in \R^d\label{32.4}
\end{equation}
where $C>0$ is the normalizing constant:
$$
C=\int_{\R^d}\prod_{j=1}^d\Big({2\sin x_k\over x_k}\Big)^2dx_1\cdots dx_d.
$$
Clearly, $h$ is symmetric. One can verify that the Fourier transform
$\widehat{h}$ is
$$
\widehat{h}(\lambda)=\int_{\R^d}h(x)e^{i\lambda\cdot x}dx
=C^{-1}(2\pi)^d\Big(1_{[-1,1]^d}\ast 1_{[-1,1]^d}\Big)(\lambda).
$$
In particular, $\widehat{h}$ is non-negative, continuous, with compact
support in the set $[-2,2]^d$,  and
\begin{equation}
\widehat{h}(\lambda)\leq \widehat{h}(0)=1.\label{p3.1}
\end{equation}
For each $\epsilon >0$, write
\begin{equation}
h_\epsilon (x)=\epsilon^{-d}h(\epsilon^{-1}x),\hskip.2in x\in \R^d.\label{ref-d}
\end{equation}

Let
\begin{equation}
\wp_{\al, \ep} ( \lambda)={C_{d,\si}\widehat{h}^{2}(\ep \lambda) 
\over \al+|\la|^{d-\si}},\label{64.2}
\end{equation}
and note that by (\ref{p3.1}) 
\begin{equation}
 \wp_{\al, \ep}( \lambda)\leq  
  \wp_{0, \ep}( \lambda)
\leq \phi_{d-\si}( \lambda).\label{64.2a}
\end{equation}
Set 
\begin{equation}
\th_{\alpha, \ep}( x)=\int_{\R^d} e^{ ix\cdot \la}\wp_{\al, \ep} ( \lambda) \,d\la,\label{64.4}
\end{equation}
and for any $A\subset R_{+}^{2}$ define
\begin{equation}
\eta_{\alpha,\ep}(A)
=\int \int_{A}
\th_{\alpha, \ep}(X_{s_1}- X_{s_2})\,ds_1\, ds_2.\label{a4.2}
\end{equation} 
Then
\begin{eqnarray}
\eta_{\alpha,  \ep}\big([0,\tau]^{2}\big)&=&\int_0^{\tau} \int_0^{\tau}\int_{\R^d} e^{ i(X_{s_1}- X_{s_2})\cdot \la}\wp_{\al, \ep} ( \lambda)\,d\la\,ds_1\, ds_2
\label{p3.11}\\
&=&\int_{\R^d}\wp_{\al, \ep} ( \lambda)\bigg\vert 
\int_0^\tau e^{i\lambda\cdot X_s}ds\bigg\vert^2d\lambda.   \nonumber
\end{eqnarray}

Let $f(\lambda)$ be a continuous and strictly positive
 function on $\R^d$ such that
$$
f(-\lambda)=f(\lambda)\hskip.1in\hbox{and}\hskip.1in
\int_{\R^d}f^2(\lambda)\wp_{\al, \ep} ( \lambda)
d\lambda =1.
$$
By the Cauchy-Schwartz inequality
$$
\eta_{\alpha,  \ep} \big([0,\tau]^2\big)^{1/2}
\ge \bigg\vert\int_{\R^d}\wp_{\al, \ep} ( \lambda) f(\lambda)
\bigg[\int_0^\tau e^{i\lambda\cdot X_s}ds\bigg]
d\lambda\bigg\vert.
$$
Hence,
$$
\begin{aligned}
&\E\Big[\eta_{\alpha,  \ep}\big([0,\tau]^2\big)^{n/2}\Big]\\
&\ge\int_{\R^n}d\lambda_1\cdots d\lambda_n\Big(\prod_{k=1}^n\wp_{\al, \ep}(\lambda_k)
f(\lambda_k)\Big)
\E\int_{[0,\tau]^n}ds_1\cdots ds_n\exp\bigg\{i\sum_{k=1}^n\lambda_k\cdot X_{s_k}
\bigg\}.
\end{aligned}
$$
Let $\Sigma_n$ be the permutation group on the set $\{1,\cdots, n\}$ and adopt
the notation
$$
[0, t]_<^n=\Big\{(s_1,\cdots, s_n)\in [0, t]^n;\hskip.1in s_1<\cdots<t_n\Big\}.
$$
We have
$$
\begin{aligned}
&\E\int_{[0,\tau]^n}ds_1\cdots ds_n\exp\bigg\{i\sum_{k=1}^n\lambda_k\cdot X_{s_k}
\bigg\}\\
&=\int_0^\infty dt e^{-t}
\int_{[0,t]^n}ds_1\cdots ds_n\E\exp\bigg\{i\sum_{k=1}^n\lambda_k\cdot X_{s_k}
\bigg\}\\
&=\sum_{\sigma\in\Sigma_n}\int_0^\infty dt e^{-t}
\int_{[0,t]_<^n}ds_1\cdots ds_n
\E\exp\bigg\{i\sum_{k=1}^n\lambda_{\sigma(k)}\cdot X_{s_k}\bigg\}.
\end{aligned}
$$
Write
$$
\sum_{k=1}^n\lambda_{\sigma(k)}\cdot X_{s_k}
=\sum_{k=1}^n\Big(\sum_{j=k}^n\lambda_{\sigma(j)}\Big)\big(X_{s_k}-X_{s_{k-1}}\big).
$$
By  independence,
$$
\E\exp\bigg\{i\sum_{k=1}^n\lambda_{\sigma(k)}\cdot X_{s_k}\bigg\}
=\exp\bigg\{-\sum_{k=1}^n(s_k-s_{k-1})\Big\vert\sum_{j=k}^n\lambda_{\sigma(j)}
\Big\vert^\beta\bigg\}.
$$
Therefore,
$$
\E\int_{[0,\tau]^n}ds_1\cdots ds_n\exp\bigg\{i\sum_{k=1}^n\lambda_k\cdot X_{s_k}
\bigg\}=\sum_{\sigma\in\Sigma_n}\prod_{k=1}^nQ\Big(\sum_{j=k}^n\lambda_{\sigma(j)}
\Big)
$$
where  $Q(\la)=\big(1+\psi(\lambda)\big)^{-1}$.

Hence,
\begin{align}\label{3.7}
&\E\Big[\eta_{\alpha, \ep}\big([0,\tau]^2\big)^{n/2}\Big]\\
&\ge\sum_{\sigma\in\Sigma_n}
\int_{\R^n}d\lambda_1\cdots d\lambda_n\Big(\prod_{k=1}^n\wp_{\al, \ep}(\lambda_k)
f(\lambda_k\Big)\prod_{k=1}^nQ\Big(\sum_{j=k}^n\lambda_{\sigma(j)}
\Big)\nn\\
&=n!\int_{\R^n}d\lambda_1\cdots d\lambda_n\Big(\prod_{k=1}^n\wp_{\al, \ep}(\lambda_k)
f(\lambda_k\Big)\prod_{k=1}^nQ\Big(\sum_{j=k}^n\lambda_{j}
\Big).\nn
\end{align}

By a change of variables 
$$
\begin{aligned}
&\int_{\R^n}d\lambda_1\cdots d\lambda_n\Big(\prod_{k=1}^n\wp_{\al, \ep}(\lambda_k)
f(\lambda_k\Big)\prod_{k=1}^nQ\Big(\sum_{j=k}^n\lambda_{j}\Big)\\
&=\int_{\R^n}d\lambda_1\cdots d\lambda_n\Big(\prod_{k=1}^n
\wp_{\al, \ep}(\lambda_k-\lambda_{k-1})
f(\lambda_k-\lambda_{k-1})Q(\lambda_k)\Big),
\end{aligned}
$$
where we follow the convention $\lambda_0=0$.

Applying an argument based on the spectral representation of self-adjoint
operators in $L^2$  (see, (3.7)-(3.10) in \cite{BCR3}) to the right 
hand side, 
\begin{align}\label{3.8}
&\liminf_{n\to\infty}{1\over n}\log{1\over n!}
\int_{\R^n}d\lambda_1\cdots d\lambda_n\Big(\prod_{k=1}^n
\wp_{\al, \ep}(\lambda_k-\lambda_{k-1})
f(\lambda_k-\lambda_{k-1})Q(\lambda_k)\Big)\nn\\
&\ge\log\sup_{\|g\|_{2}=1}\int_{\R^d}\wp_{\al, \ep}(\lambda)f(\lambda)\bigg[
\int_{\R^d}\sqrt{Q(\lambda+\gamma)}
\sqrt{Q(\gamma)}g(\lambda+\gamma)g(\gamma)d\gamma\bigg]d\lambda.
\end{align}

Write
\begin{align}\label{3.1}
\rho_{\alpha,\ep}=\sup_{\vert\vert g\vert\vert_2=1}
\int_{\R^d}\wp_{\al, \ep}( \lambda)
\bigg[\int_{\R^d}\sqrt{Q(\lambda+\gamma)}
\sqrt{Q(\gamma)}g(\lambda+\gamma)g(\gamma)d\gamma\bigg]^2d\lambda.
\end{align}
Taking the supremum over $f$   on the right hand side of (\ref{3.8})
leads to the conclusion that
\begin{equation}
\liminf_{n\to\infty}{1\over n}\log{1\over n!}
\E\Big[\eta_{\alpha, \ep}\big([0,\tau]^2\big)^{n/2}\Big]\geq {1\over 2}\log\rho_{\alpha, \ep}.\label{p3.12}
\end{equation}
\qed

We note for future reference that it follows easily from 
(\ref{p3.1}) and  \cite[Lemma 1.6]{BCR3} that
\bea
&&
0\leq \rho-\rho_{\alpha, \ep}\leq \sup_{\vert\vert g\vert\vert_2=1}
\int_{\R^d}\(1-{\widehat{h}^{2}(\ep \lambda)|\la|^{d-\si} 
\over \al+|\la|^{d-\si}}\)\phi_{d-\si}( \lambda)\label{}\\
&&\hspace{1.5 in}
\bigg[\int_{\R^d}\sqrt{Q(\lambda+\gamma)}
\sqrt{Q(\gamma)}g(\lambda+\gamma)g(\gamma)d\gamma\bigg]^2d\lambda\nonumber\\
&&\leq C\sup_{\vert\vert g\vert\vert_2=1}\Bigg\vert\Bigg\vert \(1-{\widehat{h}^{2}(\ep \lambda)|\la|^{d-\si} 
\over \al+|\la|^{d-\si}}\)g\Bigg\vert\Bigg\vert_2\,\,,\nn
\eea
so that 
\begin{equation}
\lim_{\alpha,\ep\rar 0}\rho_{\alpha, \ep}=\rho.\label{p3.0}
\end{equation}

Note also  for future reference that for any $c>0$ and $t>0$,
\begin{eqnarray}
\eta_{\alpha, \ep}\big([0,ct]_<^2\big)&=&\int_0^{ct} \int_0^{s_2}\int_{\R^d} e^{ i(X_{s_1}- X_{s_2})\cdot \la}{C_{d,\si}\widehat{h}^{2}(\ep \lambda) 
\over \alpha+|\la|^{d-\si}}\,d\la\,ds_1\, ds_2
\label{p3.2}\\
&=& c^{2}  \int_0^{t} \int_0^{s_2}\int_{\R^d} e^{ i(X_{cs_1}- X_{cs_2})\cdot \la}{C_{d,\si}\widehat{h}^{2}(\ep \lambda) 
\over \al+|\la|^{d-\si}}\,d\la\,ds_1\, ds_2\nonumber\\
&\buildrel d\over =& c^{2}  \int_0^{t} \int_0^{s_2}\int_{\R^d} e^{ i(X_{s_1}- X_{s_2})\cdot c^{1/\bb}\la}{C_{d,\si}\widehat{h}^{2}(\ep \lambda) 
\over \alpha+|\la|^{d-\si}}\,d\la\,ds_1\, ds_2\nonumber\\
&  =& c^{2-d/ \beta}  \int_0^{t} \int_0^{s_2}\int_{\R^d} e^{ i(X_{s_1}- X_{s_2})\cdot \la}{C_{d,\si}\widehat{h}^{2}(\ep \lambda/c^{1/\bb}) 
\over \alpha+|\la/c^{1/\bb}|^{d-\si}}\,d\la\,ds_1\, ds_2\nonumber\\
&  =& c^{2  -\sigma/ \beta}   \int_0^{t} \int_0^{s_2}\int_{\R^d} e^{ i(X_{s_1}- X_{s_2})\cdot  \la}{C_{d,\si}\widehat{h}^{2}(\ep \lambda/c^{1/\bb}) 
\over \alpha c^{(d-\si)/\bb}+|\la|^{d-\si}}\,d\la\,ds_1\, ds_2\nonumber\\
&=&c^{2  -\sigma/ \beta}
\eta_{\alpha c^{(d-\si)/\bb},\,\epsilon c^{-1/\beta}}\big([0, t]_<^2\big).\nn
\end{eqnarray}

\begin{lemma}\label{3.9}
\begin{align}\label{32.1}
\limsup_{n\to\infty}{1\over n}\log {1\over n!}
\E\Big[\eta_{\alpha, \ep}\big([0,\tau]^2\big)^{n/2}\Big]
\le {1\over 2}\log\rho_{\alpha, \ep}.
\end{align}
\end{lemma}

\Proof 
Let
\begin{equation}
\bar\th_{ \al,\ep}( x)
=\int_{\R^d} e^{ ix\cdot \la}{C_{d,\si}\widehat{h}(\ep \lambda) 
\over \al+|\la|^{d-\si}}  \,d\la,\label{64.3}
\end{equation}
and
\begin{equation}
\psi^{x}_{ \al,\ep}\big(A\big)
=\int \int_{A}
\bar\th_{ \al,\ep}(X_{s_1}- X_{s_2}-x)\,ds_1\, ds_2.\label{a4.2}
\end{equation} 
Then
\begin{eqnarray}
\psi^{x}_{ \al,\ep}\big([0,\tau]^{2}\big)&=&\int_0^{\tau} \int_0^{\tau}\int_{\R^d} e^{ i(X_{s_1}- X_{s_2}-x)\cdot \la}{C_{d,\si}\widehat{h}(\ep \lambda) 
\over \al+|\la|^{d-\si}}\,d\la\,ds_1\, ds_2
\nn\\
&=&\int_{\R^d}e^{-ix\cdot \la} {C_{d,\si}\widehat{h}(\ep \lambda) 
\over \al+|\la|^{d-\si}}  \bigg\vert 
\int_0^\tau e^{i\lambda\cdot X_s}ds\bigg\vert^2d\lambda,   \label{p3.11}
\end{eqnarray}
so that
\begin{equation}
\eta_{ \al,\ep}\big([0,\tau]^{2}\big)
=\int h_{\ep}(x)\psi^{x}_{ \al,\ep }\big([0,\tau]^{2}\big)\,dx.\label{32.11m1}
\end{equation}
Let $M>0$ be a large but fixed number.
Define the random measure $\tilde{\psi}_{ \al,\ep }^x(\cdot)$ as
$$
\tilde{\psi}_{ \al,\ep }^x(A)=\sum_{y\in\Z^d}\psi_{ \al,\ep }^{yM+x}(A)
$$
and write
$$
\tilde{h}_{\ep}(x)=\sum_{y\in\Z^d}h_{\ep} (yM+x).
$$
Then
$$
\begin{aligned}
\eta_{ \al,\ep }([0, \tau]^2\big)
&=\sum_{y\in\Z^d}\int_{[0, M]^2}h_{\ep}(yM+x)\psi_{ \al,\ep }^{yM+x}
\big([0, \tau]^2\big)\\
&\le \int_{[0, M]^d}\tilde{h}_{\ep}(z)
\tilde{\psi}_{ \al,\ep }^z([0, \tau]^2\big)dz.
\end{aligned}
$$

In addition, by Parseval identity,
$$
\begin{aligned}
&\int_{[0, M]^d}\tilde{h}_{\ep}(z)
\tilde{\psi}_{ \al,\ep }^z([0, \tau]^2\big)dz\\
&=M^{-d}\sum_{y\in\Z^d}\bigg(\int_{[0, M]^d}\tilde{h}_{\ep}(x)
\exp\Big\{-i{2\pi\over M}(y\cdot x)\Big\}dx\bigg)\\
&\times\bigg(\int_{[0, M]^d}\tilde{\psi}_{ \al,\ep }^x ([0, \tau]^2\big)
\exp\Big\{i{2\pi\over M}(y\cdot x)\Big\}dx\bigg).
\end{aligned}
$$
By the periodicity of $\tilde{h}_{\ep}$,
$$
\begin{aligned}
&\int_{[0, M]^d}\tilde{h}_{\ep}(x)
\exp\Big\{-i{2\pi\over M}(y\cdot x)\Big\}dx\\
&=\sum_{z\in\Z^d}\int_{[0, M]^d}h_{\ep}(zM+x)
\exp\Big\{-i{2\pi\over M}(y\cdot x)\Big\}dx\\
&=\sum_{z\in\Z^d}\int_{[0, M]^d}h_{\ep'}(zM+x)
\exp\Big\{-i{2\pi\over M}\big(y\cdot (zM+x)\big)\Big\}dx\\
&=\int_{\R^d}h_{\ep}(x)
\exp\Big\{-i{2\pi\over M}(y\cdot x)\Big\}dx
=\hat{h}\Big({2\pi\epsilon\over M}y\Big).
\end{aligned}
$$
Similarly,
$$
\begin{aligned}
&\int_{[0, M]^d}\tilde{\psi}^x_{ \al,\ep }([0, \tau]^2\big)
\exp\Big\{i{2\pi\over M}(y\cdot x)\Big\}dx\\
&=\int_{\R^d}\psi^x_{ \al,\ep }([0, \tau]^{2}\big)\exp\Big\{-i{2\pi\over M}(y\cdot x)\Big\}dx\\
&=(2\pi)^d{C_{d,\si}\widehat{h}\Big({2\pi\epsilon\over M}y\Big) 
\over \al+|{2\pi\over M}y|^{d-\si}}\bigg\vert\int_0^\tau 
\exp\Big\{i{2\pi\over M}\big(y\cdot X_s\big)\Big\}ds\bigg\vert^2,
\end{aligned}
$$
where the last step follows by (\ref{p3.11}).

Summarizing the computation,
\begin{align}\label{3.11}
&\eta_{ \al,\ep }([0, \tau]^2\big)\\
&\le \Big({2\pi\over M}\Big)^d
\sum_{y\in\Z^d} \wp_{\al, \ep}\Big({2\pi\over M}y\Big)\bigg\vert\int_0^\tau 
\exp\Big\{i{2\pi\over M}(y\cdot X_s)\Big\}ds\bigg\vert^2\nn\\
&= \Big({2\pi\over M}\Big)^d\sum_{y\in E}
\wp_{\al, \ep}\Big({2\pi\over M}y\Big)\bigg\vert\int_0^\tau 
\exp\Big\{i{2\pi\over M}(y\cdot X_s)\Big\}ds\bigg\vert^2,\nn
\end{align}
where 
$$
E=\Z^d\cap \Big[-{M\over\pi\epsilon}, \hskip.05in {M\over\pi\epsilon}\Big],
$$
and the last step follows from the fact that $\hat{h}(\lambda)$ is supported
on $[-2, 2]^d$.

Write 
$$
\pi_{\alpha,\ep ,M}(y)= 
\wp_{\al, \ep}\Big({2\pi\over M}y\Big).
$$
Let ${\cal H}$ be the finite dimensional $U$-space of the complex-valued
functions $g(y)$ on $E$ with 
$$
\|g\|=\bigg\{\sum_{y\in E}\vert g(y)\vert^2\pi_{\alpha,\ep ,M}(y)\bigg\}^{1/2}.
$$
Let the subset ${\cal U}\subset {\cal H}$ be defined by the property
$$
\overline{g(y)}=g(-y)\hskip.2in y\in E.
$$
Let $\delta>0$ be a small but fixed number. There are 
$f_1,\cdots, f_l\in{\cal U}$ such that
$$
\| g\|\le (1+\delta)\max_{1\le j\le l}\vert\langle f_j, g\rangle\vert
\hskip.2in g\in {\cal U}.
$$

In particular,
$$
\begin{aligned}
&\bigg\{\sum_{y\in E}\pi_{\alpha,\ep ,M}(y)\bigg\vert\int_0^\tau 
\exp\Big\{i{2\pi\over M}(y\cdot X_s)\Big\}ds\bigg\vert^2\bigg\}^{1/2}\\
&\le (1+\delta)\max_{1\le j\le l}\bigg\vert\sum_{y\in E}\pi_{\alpha,\ep ,M}(y)f_j(y)
\int_0^\tau 
\exp\Big\{i{2\pi\over M}(y\cdot X_s)\Big\}ds\bigg\vert\\
&= (1+\delta)\max_{1\le j\le l}\bigg\vert\sum_{y\in \Z^d}\pi_{\alpha,\ep ,M}(y)f_j(y)
\int_0^\tau 
\exp\Big\{i{2\pi\over M}(y\cdot X_s)\Big\}ds\bigg\vert.
\end{aligned}
$$
By (3.11),
\begin{align}\label{3.12}
&\E\Big[\eta_{\alpha,\ep }\big([0,\tau]^2\big)^{n/2}\Big]\\
&\le (1+\delta)^{n}\Big({2\pi\over M}\Big)^{nd/2}\sum_{j=1}^l
\E\bigg\vert\sum_{y\in \Z^d}\pi_{\alpha,\ep ,M}(y)f_j(y)
\int_0^\tau 
\exp\Big\{i{2\pi\over M}(y\cdot X_s)\Big\}ds\bigg\vert^n.\nn
\end{align}

We now intend to establish
\begin{align}\label{3.13}
&\limsup_{n\to\infty}{1\over n}\log {1\over n!}
\E\bigg\vert\sum_{y\in \Z^d}\pi_{\alpha,\ep ,M}(y)f_j(y)
\int_0^\tau 
\exp\Big\{i{2\pi\over M}(y\cdot X_s)\Big\}ds\bigg\vert^n\\
&\hspace{3.5 in}\le{1\over 2}\log \rho_{\alpha,\ep ,M}\nn
\end{align}
for each $j=1,\cdots, l$, where
\begin{align}\label{3.14}
&\rho_{\alpha,\ep ,M}=\sup_{\vert g\vert_2=1}
\sum_{x\in \Z^d}\wp_{\al, \ep}\Big({2\pi\over M}x\Big)\\
&\hspace{1 in}\times\bigg[\sum_{y\in\Z^d}\sqrt{Q\Big({2\pi\over M}(x+y)\Big)}
\sqrt{Q\Big({2\pi\over M}y\Big)}g(x+y)g(y)\bigg]^2,\nn
\end{align}
  where the supremum is   over all functions $g(x)$ on
$\Z^d$ satisfying
$$
\vert g\vert_2\equiv\bigg\{\sum_{y\in \Z^d}g^2(y)\bigg\}^{1/2}=1.
$$

Write $f=f_j$. Using  estimates of the form
$$
\E \vert X\vert^{2m+1}\le\Big\{\E \vert X\vert^{2m}\Big\}^{1/2}
\Big\{\E \vert X\vert^{2(m+1)}\Big\}^{1/2},
$$
we need only  consider the case $n=2m$. That is, we need only to show that
\begin{align}
&\limsup_{m\to\infty}{1\over 2m}\log {1\over (2m)!}
\E\bigg\vert\sum_{y\in \Z^d}\pi_{\alpha,\ep ,M}(y)f_j(y)
\int_0^\tau 
\exp\Big\{i{2\pi\over M}(y\cdot X_s)\Big\}ds\bigg\vert^{2m}\nn
\\
&\hspace{3.3 in}\le {1\over 2}\log\rho_{\alpha,\ep ,M}.\label{3.15}
\end{align}

Indeed,
$$
\begin{aligned}
&\E\bigg\vert\sum_{y\in \Z^d}\pi_{\alpha,\ep ,M}(x)f(y)\int_0^\tau 
\exp\Big\{i{2\pi\over M}(y\cdot X_s)\Big\}ds\bigg\vert^{2m}\\
&=\sum_{y_1,\cdots, y_{2m}\in \Z^d}\bigg(
\prod_{k=1}^m\pi_{\alpha,\ep ,M}(y_k)\pi_{\alpha,\ep ,M}(y_{m+k})
f(y_k)\overline{f(y_{m+k})}\bigg)\\
&\times\E\int_{[0,\tau]^{2m}}ds_1\cdots ds_{2m}
\exp\bigg\{i{2\pi\over M}\sum_{k=1}^m\Big[
(y_k\cdot X_{s_k})-(y_{m+k}\cdot X_{s_{m+k}})\Big]\bigg\}\\
&=\sum_{y_1,\cdots, y_{2m}\in \Z^d}\bigg(
\prod_{k=1}^{2m}\pi_{\alpha,\ep ,M}(y_k)f(y_k)\bigg)\\
&\hspace{1.6 in}\E\int_{[0,\tau]^{2m}}
ds_1\cdots ds_{2m}
\exp\bigg\{i{2\pi\over M}\sum_{k=1}^{2m}
(y_k\cdot X_{s_k})\bigg\}.
\end{aligned}
$$

Similar to the computation for (\ref{3.8}), (with $n=2m$), 
the right hand 
side is equal to 
$$
(2m)!
\sum_{y_1,\cdots, y_{2m}\in \Z^d}\Big(\prod_{k=1}^{2m}
\pi_{\alpha,\ep ,M}(y_k-y_{k-1})f(y_k-y_{k-1})Q(y_k)\Big).
$$

Observe that  the same
argument of spectral representation used in the proof of (\ref{3.8}) gives
$$
\begin{aligned}
&\lim_{m\to\infty}
{1\over 2m}\log
\sum_{y_1,\cdots, y_{2m}\in \Z^d}\Big(\prod_{k=1}^{2m}\pi_{\alpha,\ep ,M}(y_k-y_{k-1})
f(y_k-y_{k-1})Q(y_k)\Big)\\&\hspace{3 in}=\log\rho_{\alpha,\ep ,M}(f)
\end{aligned}
$$
where 
$$
\begin{aligned}
&\rho_{\alpha,\ep ,M}(f)\equiv\sup_{\vert g\vert_2=1}
\sum_{x, y\in \Z^d}\pi_{\alpha,\ep ,M}(x-y)f(x-y)\sqrt{Q\Big({2\pi\over M}x\Big)}
\sqrt{Q\Big({2\pi\over M}y\Big)}g(x)g(y)\\
&= \sup_{\vert g\vert_2=1}
\sum_{x\in \Z^d}\pi_{\alpha,\ep ,M}(x)f(x)\bigg[\sum_{y\in\Z^d}\sqrt{Q\Big({2\pi\over M}(x+y)\Big)}
\sqrt{Q\Big({2\pi\over M}y\Big)}g(x+y)g(y)\bigg].
\end{aligned}
$$
 Hence, (\ref{3.15}) follows from the fact that
$\rho_{\alpha,\ep ,M}(f)\le\sqrt{\rho_{\alpha,\ep ,M}}$, (Cauchy-Schwartz).

By (\ref{3.12}) and (\ref{3.13}),
$$
\limsup_{n\to\infty}{1\over n}\log {1\over n!}
\E\Big[\eta_{\alpha,\ep }\big([0,\tau]^2\big)^{n/2}\Big]
\le {1\over 2}\log\bigg\{\Big({2\pi\over M}\Big)^d\rho_{\alpha,\ep ,M}\bigg\}.
$$

We show in Theorem \ref{theo-alt} below that
$$
\limsup_{M\to\infty}\Big({2\pi\over M}\Big)^d\rho_{\alpha,\ep,M}\le\rho_{\alpha,\ep }.
$$
This will then show that
$$
\limsup_{n\to\infty}{1\over n}\log {1\over n!}
\E\Big[\eta_{\alpha,\ep }\big([0,\tau]^2\big)^{n/2}\Big]
\le {1\over 2}\log\rho_{\alpha,\ep }.
$$
\qed

Combining (\ref{p3.12}) and Lemma \ref{3.9}, we have shown that
\begin{align}\label{3.16}
\lim_{n\to\infty}{1\over n}\log {1\over n!}
\E\Big[\eta_{\alpha,\ep }\big([0,\tau]^2\big)^{n/2}\Big]
= {1\over 2}\log\rho_{\alpha,\ep }.
\end{align}

By Taylor's expansion,
\begin{align}\label{3.17}
\E\exp\bigg [\theta\Big(\eta_{\alpha,\ep }\big([0,\tau]^2\big)\Big)^{1/2}\bigg ]
\left\{\begin{array}{ll}<\infty\hskip.2in\mbox{for}\hskip.2in \theta<\rho_{\al, \ep}^{-1/2}\\\\
=\infty\hskip.2in\mbox{for}\hskip.2in\theta>\rho_{\al, \ep}^{-1/2}.\end{array}\right.
\end{align}

In addition, replacing $\tau$ by $t$ in (\ref{p3.11}),
$$
\eta_{\al, \ep} \big([0, t]^2\big)
=\int_{\R^d}\wp_{\al, \ep}(\lambda)\bigg\vert 
\int_0^t e^{i\lambda\cdot X_u}du\bigg\vert^2d\lambda.
$$

By the  triangle inequality one has that for any $s, t>0$
$$
\begin{aligned}
&\bigg\{\int_{\R^d}\wp_{\al, \ep}(\lambda)\bigg\vert 
\int_0^t e^{i\lambda\cdot X_u}du\bigg\vert^2d\lambda\bigg\}^{1/2}\\
&\le\bigg\{\int_{\R^d}\wp_{\al, \ep}(\lambda)\bigg\vert 
\int_0^s e^{i\lambda\cdot X_u}du\bigg\vert^2d\lambda\bigg\}^{1/2}\\
&+\bigg\{\int_{\R^d}\wp_{\al, \ep}(\lambda)\bigg\vert 
\int_s^{s+t}e^{i\lambda\cdot X_u}du\bigg\vert^2d\lambda\bigg\}^{1/2}.
\end{aligned}
$$
Notice that the random quantity
$$
\begin{aligned}
&\int_{\R^d}\wp_{\al, \ep}(\lambda)\bigg\vert 
\int_s^{s+t}e^{i\lambda\cdot X_u}du\bigg\vert^2d\lambda\\
&=\int_{\R^d}\wp_{\al, \ep}(\lambda)\bigg\vert 
\int_0^{t}e^{i\lambda\cdot (X_{s+u}-X_s)}du\bigg\vert^2d\lambda
\end{aligned}
$$
is independent of $\eta_{\alpha,\ep }\big([0, s]^2\big)$ and has the
same distribution as $\eta_{\alpha,\ep }\big([0, t]^2\big)$.
Consequently, for any $\theta>0$,
$$
\begin{aligned}
&\log \E\exp\bigg\{\theta\Big(\eta_{\alpha,\ep }\big([0,s+t]^2\big)\Big)^{1/2}\bigg\}\\
&\le\log \E\exp\bigg\{\theta\Big(\eta_{\alpha,\ep }\big([0,s]^2\big)\Big)^{1/2}\bigg\}
+\log \E\exp\bigg\{\theta\Big(\eta_{\alpha,\ep }\big([0,t]^2\big)\Big)^{1/2}\bigg\}.
\end{aligned}
$$
Thus, the limit
\be
L(\th):=\lim_{t\to\infty}
{1\over t}
\log \E\exp\bigg\{\theta\Big(\eta_{\alpha,\ep }\big([0,t]^2\big)\Big)^{1/2}\bigg\}\label{p3.17}
\ee
exists as extended real number. In view of the relation
$$
\E\exp\bigg\{\theta\Big(\eta_{\alpha,\ep }\big([0,\tau]^2\big)\Big)^{1/2}\bigg\}
=\int_0^\infty dt e^{-t}
\E\exp\bigg\{\theta\Big(\eta_{\alpha,\ep }\big([0,t]^2\big)\Big)^{1/2}\bigg\},
$$
and by (\ref{3.17}) and   the relation $\eta_{\alpha,\ep }\big([0,t]^2\big)
=2\eta_{\alpha,\ep }\big([0,t]^2_<\big)$,
\begin{align}\label{3.18}
\lim_{t\to\infty}{1\over t}
\log \E\exp\bigg\{\th
\Big(\eta_{\alpha,\ep }\big([0,t]_<^2\big)\Big)^{1/2}\bigg\}\left\{\begin{array}{ll}\leq 1\hskip.2in\mbox{for}\hskip.2in \theta<\sqrt{2\over\rho_{\alpha,\ep }}\\\\
\geq 1\hskip.2in\mbox{for}\hskip.2in\theta>\sqrt{2\over\rho_{\alpha,\ep }}.\end{array}\right.
\end{align} 

Let $(\widetilde{X}_t, \tilde{\tau})$ be an independent copy of
$(X_t, \tau)$
and recall that the random measure $\zeta(\cdot)$ is defined in (\ref{1.4}).
Define
$$
\zeta^x(A)=\int\!\!\int_A\vert X_s- \widetilde{X}_t-x\vert^{-\sigma}dsdt,
\hskip.2in A\subset (\R^+)^2
$$
\begin{equation}\label{3.19}
\zeta_{\alpha,\ep }(A)=\int_{\R^d}\th_{\alpha,\ep } (x)\zeta^x(A)dx.
\end{equation} 

By \cite[Lemma 5.1]{BCR3}, (with $p=2$),
$$
\lim_{\al,\epsilon\to 0^+}\limsup_{n\to\infty}{1\over n}\log{1\over (n!)^2}
\E\Big[\zeta\big([0,\tau]\times[0,\tilde{\tau}]\big)
-\zeta_{\alpha,\ep }\big([0,\tau]\times[0,\tilde{\tau}]\big)\Big]^n=-\infty.
$$
This leads, by an similar argument, to the fact   that for any $\theta>0$
\begin{equation}\label{3.20}
\lim_{\al,\epsilon\to 0^+}\limsup_{t\to\infty}
{1\over t}\log\E\exp\bigg\{\theta\Big\vert\zeta\big([0,t]^2\big)
-\zeta_{\alpha,\ep }\big([0,t]^2\big)\Big\vert^{1/2}\bigg\}=0.
\end{equation}

\section{Large deviations}\label{sec-largedev}

In spite of their similarity, the large deviations in
(\ref{2.2}) and (\ref{2.3}) require different strategies.
When $0<\sigma<\min\{\beta, d\}$, both parts, off and near the time
diagonal, make contribution to the large deviation
given in (\ref{2.2}). When 
$\displaystyle \beta\le\sigma<\min\Big\{{3\over 2}\beta, d\Big\}$,
on the other hand, renormalization makes the off-diagonal part
the only source that contributes
to the large deviation
given in (\ref{2.3}). Accordingly, different proofs are given
for the two  cases.

\subsection{Proof of Theorem \ref{theo-ldp-1}, $0<\sigma <\min\{\beta, d\}$}

Recall (\ref{eta-3})-(\ref{eta-4}). 
Since 
\begin{align}\label{eta-4x}
\eta^{z}\big([0, t]^2\big)=C\int_{\R^d}\xi(t, x) \xi(t, x+z)  dx, 
\end{align}
and translation is continuous in $L^{2}(\R^d)$, we see that $\eta^{z}\big([0, t]^2\big)
$ is continuous in $z$, almost surely. Hence if  $f(x)\in \mathcal{S}(R^{d})$   with $\int f(x)\,dx=1$, and defining $f_{\de}$ as in (\ref{ref-d}) we have
\begin{equation}
\lim_{\de \rar 0}\int   \( \int\!\!\int_{[0, t]^2}\vert X_s-X_r-x\vert^{-\sigma}\,drds\) f_{\de}(x)\,dx=\eta\big([0, t]^2\big).\label{p5.2}
\end{equation}
But
\begin{eqnarray}
&&\int   \( \int\!\!\int_{[0, t]^2}\vert X_s-X_r-x\vert^{-\sigma}\,drds\) f_{\de}(x)\,dx
\label{p5.3}\\
&&= \int\!\!\int_{[0, t]^2}   \(\int  \vert X_s-X_r-x\vert^{-\sigma}f_{\de}(x)\,dx\)    \,drds\nonumber\\
&&=  \int  \vert x\vert^{-\sigma}F_{\de}(x)\,dx, \nonumber
\end{eqnarray}
where
\begin{equation}
F_{\de}(x)= \int\!\!\int_{[0, t]^2}  f_{\de}(x+X_s-X_r)    \,drds \label{}
\end{equation}
is in $\mathcal{S}(R^{d})$ with 
\begin{equation}
\wh F_{\de}(x)= \int\!\!\int_{[0, t]^2}  e^{i(X_s-X_r)\cdot\la}  \wh f (\de\,\la)   \,drds
=\wh f (\de\,\la)\, 
\bigg\vert\int_0^t e^{i\lambda\cdot X_s}ds\bigg\vert^2.\label{p5.3b}
\end{equation}
Hence using (\ref{a1.4})
\begin{eqnarray}
&&  \int  \vert x\vert^{-\sigma}F_{\de}(x)\,dx =\int_{\R^d}\phi_{d-\sigma}(\lambda)\,\wh f (\de\,\la)\,
\bigg\vert\int_0^t e^{i\lambda\cdot X_s}ds\bigg\vert^2d\lambda.\label{p5.3c}
\end{eqnarray}
This shows that
\begin{equation}
\eta\big([0, t]^2\big)=\int_{\R^d} \phi_{d-\sigma}(\lambda)
\bigg\vert\int_0^t e^{i\lambda\cdot X_s}ds\bigg\vert^2d\lambda.\label{}
\end{equation}
Hence using (\ref{p3.11})
$$
\eta\big([0, t]^2\big)-\eta_{\al,\ep}\big([0, t]^2\big)
=\int_{\R^d}\big[\phi_{d-\sigma}(\lambda)-\wp_{\al,\ep}( \lambda)\big]
\bigg\vert\int_0^t e^{i\lambda\cdot X_s}ds\bigg\vert^2d\lambda.
$$
Note that   
$\eta\big([0, t]^2\big)-\eta_{\al,\ep}\big([0, t]^2\big)\ge 0$ by (\ref{64.2a}).

As in the proof of (\ref{p1.5}), $\eta\big([0, t]^2\big)-\eta_\epsilon\big([0, t]^2\big)$
is sub-additive. Hence, for any $\theta>0$ 
$$
\begin{aligned}
&\lim_{t\to\infty}{1\over t}\E\exp\bigg\{\theta\Big(\eta\big([0, t]^2\big)
-\eta_{\al,\ep}\big([0, t]^2\big)\Big)^{1/2}\bigg\}\\
&=\inf_{T>0}\log\E\exp\bigg\{\theta\bigg(\int_{\R^d}\big[\phi_{d-\sigma}(\lambda)-\wp_{\al,\ep}( \lambda)\big]
\bigg\vert\int_0^Te^{i\lambda\cdot X_s}ds\bigg\vert^2d\lambda\bigg)^{1/2}\bigg\}\\
&\le\log\E\exp\bigg\{\theta\bigg(\int_{\R^d}\big[\phi_{d-\sigma}(\lambda)-\wp_{\al,\ep}( \lambda)\big]
\bigg\vert\int_0^1e^{i\lambda\cdot X_s}ds\bigg\vert^2d\lambda\bigg)^{1/2}\bigg\}.
\end{aligned}
$$
Applying the dominated convergence theorem (based on Theorem \ref{eta-1})
to the right hand side leads to
\begin{align}\label{4.1'}
\lim_{\al, \epsilon\to 0^+}
\lim_{t\to\infty}{1\over t}\E\exp\bigg\{\theta\Big(\eta\big([0, t]^2\big)
-\eta_{\al,\ep}\big([0, t]^2\big)\Big)^{1/2}\bigg\}=0
\end{align}
for each $\theta>0$.

Using (\ref{3.18}) and (\ref{p3.0})
we obtain that
$$
\limsup_{t\to\infty}{1\over t}\log\E\exp\bigg\{\theta_{1}
\Big\vert\eta\big([0,t]_<^2\big)\Big\vert^{1/2}\bigg\}
\le 1,\hskip.2in \theta_{1}<\sqrt{2\over \rho},
$$
$$
\liminf_{t\to\infty}{1\over t}\log\E\exp\bigg\{\theta_{2}
\Big\vert\eta\big([0,t]_<^2\big)\Big\vert^{1/2}\bigg\}
\ge 1,\hskip.2in \theta_{2}>\sqrt{2\over \rho}.
$$
For any $\theta>0$,   using  the substitutions
$$
t=a^{\beta/\sigma}\Big({\theta\over\theta_{1}}\Big)^{2\over  2-\sigma/\beta}
\hskip.1in
\hbox{and}
\hskip.1in
t=a^{\beta/\sigma}\Big({\theta\over\theta_{2}}\Big)^{2\over  2-\sigma/\beta}
$$
together with the scaling (\ref{1.3})  we obtain
$$
\limsup_{a\to\infty}a^{-\beta/\sigma}\log\E\exp\bigg\{\theta 
a^{\beta/\sigma-1/2}
\Big\vert\eta\big([0,1]_<^2\big)\Big\vert^{1/2}\bigg\}
\le \Big({\theta\over\theta_{1}}\Big)^{2\over  2-\sigma/\beta},
$$
$$
\liminf_{a\to\infty}a^{-\beta/\sigma}\log\E\exp\bigg\{\theta 
a^{\beta/\sigma-1/2}
\Big\vert\eta\big([0,1]_<^2\big)\Big\vert^{1/2}\bigg\}
\ge \Big({\theta\over\theta_{2}}\Big)^{2\over  2-\sigma/\beta}.
$$
Letting $\theta_{1}, \theta_{2}\to \sqrt{2/\rho}$ gives
\begin{align}\label{4.2'}
\lim_{a\to\infty}a^{-\beta/\sigma}\log
\E\exp\bigg\{\theta a^{\beta/\sigma-1/2} \Big(
\eta\big([0,1]_<^2\big)\Big)^{1/2}
\bigg\}=\theta^{2\over  2-\sigma/\beta}(\rho/2)^{ 1 \over  2-\sigma/\beta}.
\end{align}
Therefore, the large deviation given in (\ref{2.2})
follows from  the  G\"artner-Ellis theorem  (Theorem 2.3.6, \cite{DZ}).
\qed

\subsection{Proof of Theorem \ref{theo-ldp-2},  $\displaystyle\beta
\le\sigma<\min\Big\{{3\over 2}\beta, d\Big\}$}

Notice that for $a>0$,
\begin{align}\label{4.1}
\P\Big\{\big\vert\gamma\big([0, 1]_<^2\big)\big\vert\ge a\Big\}
=\P\Big\{\gamma\big([0, 1]_<^2\big)\ge a\Big\}
+\P\Big\{-\gamma\big([0, 1]_<^2\big)\ge a\Big\}.
\end{align}
We claim that
\begin{align}\label{4.2}
\lim_{a\to\infty}a^{-\beta/\sigma}\log
\P\Big\{-\gamma\big([0, 1]_<^2\big)\ge a\Big\}
=-\infty.
\end{align}

Let $m\ge 1$ be a fixed but arbitrary integer and write
$$
D_m=\bigcup_{k=1}^{m-1}\Big[{k-1\over m}, {k\over m}\Big]
\times\Big[{k\over m}, 1\Big].
$$
We note that
$$
\begin{aligned}
\gamma\big([0, 1]_<^2\big)&=\eta(D_m)-\E\eta(D_m)+
\sum_{k=1}^m\gamma\Big(\Big[{k-1\over m}, {k\over m}
\Big]_<^2\Big)\\
&\ge
-\E\eta(D_m)+
\sum_{k=1}^m\gamma\Big(\Big[{k-1\over m}, {k\over m}
\Big]_<^2\Big)\\
&\buildrel d\over =-\E\eta(D_m)+m^{-(2 -\sigma/\beta)}
\sum_{k=1}^m\gamma\big([k-1, k]_<^2\big).
\end{aligned}
$$
By (\ref{1.7}) we see that $\E\eta(D_m)<\infty$ for $m$ fixed.
Thus,
$$
\begin{aligned}
&\limsup_{a\to\infty}a^{-\beta/\sigma}\log
\P\Big\{-\gamma\big([0, 1]_<^2\big)\ge a\Big\}\\
&\leq \limsup_{a\to\infty}a^{-\beta/\sigma}\log
\P\Big\{-\sum_{k=1}^m\gamma\big([k-1, k]_<^2\big)\ge
m^{2 -\sigma/\beta}a\Big\}.
\end{aligned}
$$
Let $c>0$ satisfy (\ref{1.12}). 
By Chebyshev's inequality,  
$$
\begin{aligned}
&\P\Big\{-\sum_{k=1}^m\gamma\big([k-1, k]_<^2\big)\ge
m^{2 -\sigma/\beta}a\Big\}\\
&\leq\P\Bigg\{ \(\sum_{k=1}^m\big\vert\gamma\big([k-1, k]_<^2\big)\big\vert\)^{\bb/\si}     \ge
m^{2 \bb/\si  -1}a^{\bb/\si} \Bigg\}\\
&\le \exp\Big\{-cm^{2 \bb/\si  -1}a^{\beta/\sigma}\Big\}\\
&\hspace{1 in}
\Bigg(\E\exp\Bigg\{c\(\sum_{k=1}^m\big\vert\gamma\big([k-1, k]_<^2\big)\big\vert\)^{\bb/\si} \Bigg\}\Bigg) \\
&\le \exp\Big\{-cm^{2\bb/\si-1}a^{\beta/\sigma}\Big\}
\bigg(\E\exp\Big\{c\big\vert\gamma\big([0,1]_<^2\big)\big\vert^{\beta/\sigma}
\Big\}\bigg)^m.
\end{aligned}
$$
For the last inequality we used   $\bb/\si<1$ and the fact that 
  the $\gamma\big([k-1, k]_<^2\big)$ are i.i.d.
Therefore,
$$
\limsup_{a\to\infty}a^{-\beta/\sigma}\log
\P\Big\{-\gamma\big([0, 1]_<^2\big)\ge a\Big\}\le 
-cm^{2\bb/\si-1}.
$$
Letting $m\to\infty$ gives (\ref{4.2}).

By (\ref{4.1}), therefore, it remains to show that
\begin{align}\label{4.3'}
\lim_{a\to\infty}a^{-\beta/\sigma}\log\P\Big\{
\big\vert\gamma\big([0,1]_<^2\big)\big\vert\ge a\Big\}
=-2^{-{\beta/\sigma}}{\sigma\over \beta}
\Big({2\beta -\sigma\over \beta}\Big)^{2\beta -\sigma\over \sigma}
\rho^{-\beta/\sigma}.
\end{align}

To this end, we  need the following lemma.

\begin{lemma}\label{4.4} There is a $C>0$ independent of $\epsilon$
 and an $\alpha >0$
such that for any $\theta>0$ 
\begin{align}\label{bound-1}
\limsup_{t\to\infty}{1\over t}\log\E\exp\bigg\{\theta\Big\vert
\gamma\big([0, t]_<^2\big)\Big\vert^{1/2}\bigg\}
\le C\theta^{2\over 2 -\sigma/\beta}
\end{align}
\begin{align}\label{bound-2}
\limsup_{t\to\infty}{1\over t}\log\E\exp\bigg\{\theta
\eta_{\al,\ep}\big([0, t]_<^2\big)^{1/2}\bigg\}
\le C\theta^{2\over 2 -\sigma/\beta}.
\end{align}
\end{lemma}

\proof Note that since $\bb/\si>2/3$, if $V\ge M^{2}t^{\sigma/\beta}$
then
\begin{eqnarray}
M^{-2(\beta/\sigma-1/2)}V^{\bb/\si}&=&M^{-2(\beta/\sigma-1/2)}V^{\bb/\si-1/2}V^{1/2}
\label{p4.1}\\
&\geq &   t^{1-\sigma/2\beta}V^{1/2}.\nonumber
\end{eqnarray}
Hence any $M>0$,
$$
\begin{aligned}
&\E\exp\Big\{\big\vert\gamma\big([0,t]_<^2\big)\big\vert^{1/2}\Big\}\\
&\le e^{  Mt}
+ \E\Bigg(
\exp\Big\{
\big\vert\gamma\big([0,t]_<^2\big)\big\vert^{1/2}\Big\};\hskip.05in
\big\vert\gamma\big([0,t]_<^2\big)\big\vert\ge M^{2}t^2\Bigg)\\
&=e^{ Mt}+
 \E\Bigg(
\exp\Big\{t^{1-\sigma/2\beta}
\big\vert\gamma\big([0,1]_<^2\big)\big\vert^{1/2}\Big\};\hskip.05in
\big\vert\gamma\big([0,1]_<^2\big)\big\vert\ge M^{2}t^{\sigma/\beta}\Bigg)\\
&\le e^{Mt}+
 \E\exp\Big\{M^{-2(\beta/\sigma-1/2)}
\big\vert\gamma\big([0,1]_<^2\big)\big\vert^{\beta/\sigma}\Big\}.
\end{aligned}
$$
By Theorem \ref{1.11}, the above estimate shows that
$$
C\equiv\limsup_{t\to\infty}{1\over t}
\log\E\exp\Big\{\big\vert\gamma\big([0,t]_<^2\big)\big\vert^{1/2}\Big\}
<\infty.
$$
Replacing $t$ by $\displaystyle\theta^{2\over 2 -\sigma/\beta}t$
yields (\ref{bound-1}).

Observe that for any $\epsilon>0$, $\rho_{\al,\ep}\le\rho$. 
Hence by (\ref{3.18})
\begin{align}\label{bound-3}
\lim_{t\to\infty}{1\over t}
\log \E\exp\bigg\{\lambda
\Big(\eta_{\al,\ep}\big([0,t]_<^2\big)\Big)^{1/2}\bigg\}\le 1
\end{align}
for any $\lambda <\sqrt{2\over\rho}$

On the other hand, by (\ref{p3.2}), for any $c>0$ and $t>0$,
\be
\eta_{\al,\ep}\big([0,ct]_<^2\big)\buildrel d\over =c^{2  -\sigma/ \beta}
\eta_{\alpha c^{(d-\si)/\bb},\,\epsilon c^{-1/\beta}}\big([0, t]_<^2\big).\label{p4.5}
\ee
Taking 
$$
c=\Big({\theta\over \lambda}\Big)^{2 \over 2 -\sigma/\beta}
$$
and replacing $t$ by $ct$, $\al$ by 
$\alpha c^{-(d-\si)/\bb}$ and $\epsilon$ by $c^{1/\beta}\epsilon$
in ({\ref{bound-3}),
$$
\lim_{t\to\infty}{1\over t}
\log \E\exp\bigg\{\theta
\Big(\eta_{\al,\ep}\big([0,t]_<^2\big)\Big)^{1/2}\bigg\}\le 
\Big({\theta\over\lambda}\Big)^{2 \over 2 -\sigma/\beta}.
$$
Letting $\lambda\to\sqrt{2\over\rho}$ leads to
$$
\lim_{t\to\infty}{1\over t}
\log \E\exp\bigg\{\theta
\Big(\eta_{\al,\ep}\big([0,t]_<^2\big)\Big)^{1/2}\bigg\}\le 
\theta^{2 \over 2 -\sigma/\beta}
 (\rho/2)^{1\over 2 -\sigma/\beta}.
$$
\qed

We now show for any $\theta>0$,
\begin{align}\label{bound-4}
\lim_{\al,\epsilon\to 0^+}
\limsup_{t\to\infty}{1\over t}\log\E\exp\bigg\{\theta\Big\vert
\gamma\big([0, t]_<^2\big)-\eta_{\al,\ep}\big([0, t]_<^2\big)\Big\vert^{1/2}
\bigg\}=0.
\end{align}

Indeed, let the integer $N\ge 1$ be large but fixed and let the sets
$$
A_l^k;\hskip.1in l=0, 1,\cdots, 2^k-1,
\hskip.1in k=0, 1,\cdots, N
$$
be defined as in (\ref{1.9}). Consider the decomposition
\begin{align}\label{4.5}
&\gamma\big([0, t]_<^2\big)-\eta_{\al,\ep}\big([0, t]_<^2\big)\\
&=\sum_{l=1}^{2^{N+1}}\big\{\gamma -\eta_{\al,\ep}\big\}
\Big(\Big[{l-1\over 2^{N+1}}t, \hskip.05in
{l\over 2^{N+1}}t\Big]_<^2\Big)
+\sum_{k=0}^N\sum_{l=1}^{2^k-1}\big\{\gamma -\eta_{\al,\ep}\big\}(A_l^k).\nn
\end{align}

Notice that for each $l=0, 1,\cdots, 2^k-1$
and $k=0, 1,\cdots, N$, and using (\ref{1.5}) and an argument similar to (\ref{p3.2})
$$
\big\{\gamma -\eta_{\al,\ep}\big\}(A_l^k)\buildrel d\over =
2^{-(k+1)(2- \sigma/\beta)}\Big\{\zeta\big([0,t]^2\big)-\E\zeta\big([0,t]^2\big)
-\zeta_{\bar\al,\bar\ep}\big([0,t]^2\big)\Big\},
$$
where $\bar\al=\alpha 2^{-(k+1)(d-\si)/\bb},\,\bar\ep=\epsilon 2^{(k+1)/\beta}$.
Also by (\ref{1.5})
$$
\E\zeta\big([0,t]^2\big)=O\Big(t^{2- \sigma/\beta}\Big).
$$
By (\ref{3.20}), using the fact that $N$ is fixed, we have that
\begin{align}\label{4.6}
\lim_{\al,\epsilon\to 0^+}
\limsup_{t\to\infty}{1\over t}\log\E\exp\bigg\{\theta\Big\vert
\sum_{k=0}^N\sum_{l=1}^{2^k-1}\big\{\gamma -\eta_{\al,\ep}\big\}(A_l^k)
\Big\vert^{1/2}\bigg\}=0.
\end{align}

Note that
$$
\gamma
\Big(\Big[{l-1\over 2^{N+1}}t, \hskip.05in
{l\over 2^{N+1}}t\Big]_<^2\Big)\buildrel d\over =
\gamma
\Big(\Big[0, \hskip.05in {t\over 2^{N+1}}\Big]_<^2\Big).
$$
Replacing $t$ by $2^{-(N+1)}t$ in (\ref{bound-1}) and (\ref{bound-2})
we have
$$
\limsup_{t\to\infty}{1\over t}
\log\E\exp\bigg\{\theta\Big\vert\gamma
\Big(\Big[{l-1\over 2^{N+1}}t, \hskip.05in
{l\over 2^{N+1}}t\Big]_<^2\Big)
\Big\vert^{1/2}\bigg\}
\le {C\theta^{2\beta\over 2\beta -\sigma}\over 2^{N+1}},
$$
$$
\limsup_{t\to\infty}{1\over t}
\log\E\exp\bigg\{\theta\Big\vert\eta_{\al,\ep}
\Big(\Big[{l-1\over 2^{N+1}}t, \hskip.05in
{l\over 2^{N+1}}t\Big]_<^2\Big)
\Big\vert^{1/2}\bigg\}
\le {C\theta^{2\beta\over 2\beta -\sigma}\over 2^{N+1}},
$$
for $l=1,\cdots 2^{N+1}$.
 
By the fact that ${2\beta\over 2\beta -\sigma}\ge 2$, and by Part (b),
Theorem 1.2.2 in \cite{Chen},
$$
\limsup_{t\to\infty}{1\over t}
\log\E\exp\bigg\{\theta\bigg(\sum_{l=1}^{2^{N+1}}\Big\vert\gamma
\Big(\Big[{l-1\over 2^{N+1}}t, \hskip.05in
{l\over 2^{N+1}}t\Big]_<^2\Big)
\Big\vert\bigg)^{1/2}\bigg\}
\le {C\theta^{2\beta\over 2\beta -\sigma}\over 2^{N+1}},
$$
$$
\limsup_{t\to\infty}{1\over t}
\log\E\exp\bigg\{\theta\bigg(\sum_{l=1}^{2^{N+1}}\Big\vert\eta_{\al,\ep}
\Big(\Big[{l-1\over 2^{N+1}}t, \hskip.05in
{l\over 2^{N+1}}t\Big]_<^2\Big)\Big\vert\bigg)^{1/2}
\bigg\}
\le {C\theta^{2\beta\over 2\beta -\sigma}\over 2^{N+1}}.
$$
Combining this with (\ref{4.5}) and (\ref{4.6}) leads to (\ref{bound-4}).

\medskip

Using (\ref{3.18}) and (\ref{bound-4}),  we find as in the proof of  (\ref{4.2'}) that for any $\th>0$
$$
\lim_{a\to\infty}a^{-\beta/\sigma}\log
\E\exp\bigg\{\theta a^{\beta/\sigma-1/2}  \Big\vert
\gamma\big([0,1]_<^2\big)\Big\vert^{1/2}
\bigg\}=\theta^{2\over  2-\sigma/\beta}(\rho/2)^{ 1 \over  2-\sigma/\beta}.
$$
Then the large deviation given in (\ref{4.3'})
follows from the G\"artner-Ellis theorem  (Theorem 2.3.6, \cite{DZ}). \qed
 
\section {Laws of the iterated logarithm}\label{sec-lil}

{\bf   Proof of Theorem \ref{theo-lil-1}:}
Using the scaling property (\ref{1.3}), our large deviation result
(\ref{2.2}) can be re-written  as
\begin{align}\label{5.1}
&\lim_{t\to\infty}(\log\log t)^{-1}\log\P\Big\{\eta\big([0, t]_<^2\big)
\ge \lambda t^{2  - \sigma/\beta}(\log\log t)^{\sigma\over\beta}\Big\}\\
&=-2^{-{\beta/\sigma}}{\sigma\over \beta}
\Big({2\beta -\sigma\over \beta}\Big)^{2\beta -\sigma\over \sigma}
\rho^{-\beta/\sigma}\lambda^{\beta/\sigma},\hskip.2in \lambda>0.\nn
\end{align}

Fix $\theta>1$  and let $t_n=\theta^n$ ($n=1,2,\cdots$).
Let the fixed numbers $\lambda_1,\lambda_2>0$ satisfy
$$
\lambda_1> 2\rho\Big({\beta\over\sigma}\Big)^{\sigma/\beta}
\Big({\beta\over 2\beta -\sigma}\Big)^{ 2\beta -\sigma\over\beta}>\lambda_2.
$$

By (\ref{5.1}),
$$
\sum_n\P\Big\{\eta\big([0, t_n]_<^2\big)
\ge \lambda_1 t_n^{2  - \sigma/\beta}(\log\log t_n)^{\sigma\over\beta}\Big\}
<\infty.
$$
Using the Borel-Cantelli lemma,
$$
\limsup_{n\to\infty}t_n^{-(2  - \sigma/\beta)}(\log\log t_n)^{-\sigma/\beta}
\eta\big([0, t_n]_<^2\big)\le\lambda_1,\hskip.2in a.s.
$$

Using the  monotonicity of $\eta\big([0, t]_<^2\big)$ we see that for any $t_n\le t\le t_{n+1}$,
$$
\begin{aligned}
&t^{-(2  - \sigma/\beta)}(\log\log t)^{-\sigma/\beta}
\eta\big([0, t]_<^2\big)\\
&\le{\displaystyle
t_{n+1}^{{2  - \sigma/\beta}}(\log\log t_{n+1})^{\sigma/\beta}\over
\displaystyle t_n^{{2  - \sigma/\beta}}(\log\log t_n)^{\sigma/\beta}
\eta\big([0, t_n]_<^2\big)}
t_{n+1}^{-(2  - \sigma/\beta)}(\log\log t_{n+1})^{-\sigma/\beta}
\eta\big([0, t_{n+1}]_<^2\big).
\end{aligned}
$$
Consequently,
$$
\limsup_{t\to\infty}t^{-(2  - \sigma/\beta)}(\log\log t)^{-\sigma/\beta}
\eta\big([0, t]_<^2\big)\le\theta^{2\beta-\sigma\over \beta}
\lambda_1,\hskip.2in a.s.
$$
Since $\theta$ can be arbitrarily close to 1 and $\lambda_1$
can be arbitrarily close to
$$
2\rho\Big({\beta\over\sigma}\Big)^{\sigma/\beta}
\Big({\beta\over 2\beta -\sigma}\Big)^{ 2\beta -\sigma\over\beta},
$$
we have proved the upper bound for (\ref{2.8}).

\medskip

On the other hand, notice that the sequence
$$
\eta\big([t_n, t_{n+1}]_<^2\big),\hskip.2in n=1,2,\cdots
$$
is independent and for each $n$,
$$
\eta\big([t_n, t_{n+1}]_<^2\big)\buildrel d\over =
\eta\big([0, t_{n+1}-t_n]_<^2\big).
$$
By (\ref{5.1}), one can make $\theta$ sufficiently large so
$$
\sum_{n}\P\Big\{\eta\big([t_n, t_{n+1}]_<^2\big)\ge\lambda_2
t_{n+1}^{2  - \sigma/\beta}(\log\log t_{n+1})^{\sigma/\beta}\Big\}=\infty.
$$
By the Borel-Cantelli lemma,
$$
\limsup_{n\to\infty}t_{n+1}^{-(2  - \sigma/\beta)}
(\log\log t_{n+1})^{-\sigma/\beta}\eta\big([t_n, t_{n+1}]_<^2\big)\ge\lambda_2,
\hskip.2in a.s.
$$
Using the  monotonicity of $\eta\big([0, t]_<^2\big)$, this leads to
$$
\limsup_{t\to\infty}t^{-(2  - \sigma/\beta)}
(\log\log t)^{-\sigma/\beta}\eta\big([0, t]_<^2\big)\ge\lambda_2,
\hskip.2in a.s.
$$
Letting
$$ 
\lambda_2\to 2\rho\Big({\beta\over\sigma}\Big)^{\sigma/\beta}
\Big({\beta\over 2\beta -\sigma}\Big)^{ 2\beta -\sigma\over\beta}
$$
yields the lower bound for (\ref{2.8}).

\medskip
{\bf   Proof of Theorem \ref{theo-lil-2}:}
 
We now turn to the proof of (\ref{2.10}). With $\lambda_1,\lambda_1>0$ as above and using (\ref{2.3}), we also have 
\begin{align}\label{5.2}
\limsup_{n\to\infty}t_n^{-(2  - \sigma/\beta)}(\log\log t_n)^{-\sigma/\beta}
\gamma\big([0, t_n]_<^2\big)\le\lambda_1,\hskip.2in a.s.
\end{align}
for any $\theta>1$; and
\begin{align}\label{5.3}
\limsup_{n\to\infty}t_{n+1}^{-(2  - \sigma/\beta)}
(\log\log t_{n+1})^{-\sigma/\beta}
\gamma\big([t_n, t_{n+1}]_<^2\big)\ge\lambda_2,\hskip.2in a.s.
\end{align}
for sufficiently large $\theta$. 

However, we can
not continue as above
due to the fact that $\gamma\big([0, t]_<^2\big)$ is not monotone in $t$. Instead, we first observe  that for any $\epsilon>0$
$$
\begin{aligned}
&\P\Big\{\sup_{t_n\le t\le t_{n+1}}\big\vert\gamma\big([0, t]_<^2\big)-
\gamma\big([0, t_n]_<^2\big)\big\vert
\ge \epsilon t_n^{2  - \sigma/\beta}(\log\log t_n)^{\sigma/\beta}\Big\}\\
&=\P\Big\{\sup_{\theta^{-1}\le t\le 1}\big\vert\gamma\big([0, t]_<^2\big)
-\gamma\big([0, \theta^{-1}]_<^2\big)\big\vert
\ge \epsilon \theta^{-{2\beta  - \sigma\over\beta}}
(\log\log t_n)^{\sigma/\beta}\Big\}.
\end{aligned}
$$
By Chebyshev's inequality and (\ref{1.14}), 
$$
\sum_n\P\Big\{\sup_{\theta^{-1}\le t\le 1}\big\vert\gamma\big([0, t]_<^2\big)
-\gamma\big([0, \theta^{-1}]_<^2\big)\big\vert
\ge \epsilon \theta^{-{2\beta  - \sigma\over\beta}}
(\log\log t_n)^{\sigma/\beta}\Big\}<\infty
$$
when $\theta$ is sufficiently close to 1. Consequently, 
\begin{align}\label{5.4}
\limsup_{n\to\infty}t_n^{-(2  - \sigma/\beta)}(\log\log t_n)^{-\sigma/\beta}
\sup_{t_n\le t\le t_{n+1}}
\big\vert\gamma\big([0, t]_<^2\big)-\gamma\big([0, t_n]_<^2\big)\big\vert\le\epsilon,\hskip.2in a.s.
\end{align}
Combining this with (\ref{5.2}) leads to the upper bound for ({2.10}).

\medskip

For the lower bound, observe that
$$
\begin{aligned}
\gamma\big([0, t_{n+1}]_<^2\big)&= \gamma\big([t_n, t_{n+1}]_<^2\big)
+\gamma\big([0, t_n]_<^2\big)+\gamma\big([0, t_n]\times[t_n, t_{n+1}]\big)\\
&\ge \gamma\big([t_n, t_{n+1}]_<^2\big)
+\gamma\big([0, t_n]_<^2\big)-\E\gamma\big([0, t_n]\times[t_n, t_{n+1}]\big).
\end{aligned}
$$
Given $\epsilon>0$, an argument similar to the one used for (\ref{5.4})
shows that
\begin{align}\label{5.5}
\limsup_{n\to\infty}t_{n+1}^{-(2  - \sigma/\beta)}
(\log\log t_{n+1})^{-\sigma/\beta}\big\vert
\gamma\big([0, t_n]_<^2\big)\big\vert\le\epsilon,\hskip.2in a.s.
\end{align}
for  $\theta$ sufficiently large.

In addition
$$
\E\gamma\big([0, t_n]\times[t_n, t_{n+1}]\big)
=\E\zeta\big([0, t_n]\times[0, t_{n+1}-t_n]\big)
=O\Big(t_{n+1}^{2  - \sigma/\beta}\Big).
$$
Combining this with (\ref{5.3}) and (\ref{5.5}) yields
$$
\limsup_{n\to\infty}t_{n+1}^{-(2  - \sigma/\beta)}
(\log\log t_{n+1})^{-\sigma/\beta}
\gamma\big([0, t_{n+1}]_<^2\big)\ge\lambda_2-\epsilon,\hskip.2in a.s. 
$$
Consequently,
$$
\limsup_{t\to\infty}t^{-(2  - \sigma/\beta)}
(\log\log t)^{-\sigma/\beta}
\gamma\big([0, t]_<^2\big)\ge\lambda_2-\epsilon,\hskip.2in a.s. 
$$
This leads to the lower bound for (\ref{2.10}). \qed

\section{The stable process in a Brownian potential}\label{sec-F}

\medskip

Throughout this section, we take $\sigma=2p-d$.

{\bf   Proof of Corollary \ref{theo-ldp-f}:}
  From (\ref{i1.8}), we have that
$$
F(1)\buildrel d\over =U\sqrt{{2\over C}\eta\big([0, 1]^2_<\big)}
$$
where $U$ is a standard normal random variable  independent of $X_t$. The remainder
of the proof for (\ref{2.6})
follows from (\ref{2.2}) and a standard computation.  \qed

{\bf   Proof of Corollary \ref{theo-lil-f}:}
Since $\{-F(t);\hskip.1in t\ge 0\}\buildrel d\over =
\{F(t);\hskip.1in t\ge 0\}$, we need only   show that
\begin{align}\label{5.6}
&\limsup_{t\to\infty}t^{-{2\beta-\sigma\over 2\beta}}
(\log\log t)^{-{\sigma+\beta\over 2\beta}} F(t)\\
&=\sqrt{8\rho\over C}\Big({\beta\over\sigma+\beta}\Big)^{\sigma+\beta\over 2\beta}
\Big({\beta\over 2\beta-\sigma}\Big)^{2\beta-\sigma\over 2\beta}\,,
\hskip.2in a.s.\nn
\end{align}

Using the scaling property (\ref{2.4}), our large deviation result
(\ref{2.6}) can be re-written  as
\begin{align}\label{f5.1}
&\lim_{t\to\infty}(\log\log t)^{-1}\log\P\Big\{F(t)
\ge \lambda t^{2\beta-\sigma\over 2\beta}(\log\log t)^{\sigma+\bb\over 2\beta}\Big\}\\
&=-{\beta+\sigma\over\beta}
\Big({C\over 8\rho}\Big)^{\beta\over \beta+\sigma}
\Big({2\beta-\sigma\over\beta}\Big)^{2\beta-\sigma\over \beta+\sigma}\lambda^{2\beta\over \sigma+\beta},\hskip.2in \lambda>0.\nn
\end{align}

Fix $\theta>1$  and let $t_n=\theta^n$ ($n=1,2,\cdots$).
Let the fixed numbers $\lambda_1,\lambda_2>0$ satisfy
$$
\lambda_1> \sqrt{8\rho\over C}\Big({\beta\over\sigma+\beta}\Big)^{\sigma+\beta\over 2\beta}
\Big({\beta\over 2\beta-\sigma}\Big)^{2\beta-\sigma\over 2\beta}>\lambda_2.
$$

By (\ref{f5.1}),
\be
\sum_n \P\Big\{F(t_{n})
\ge \lambda_{1} t_{n}^{2\beta-\sigma\over 2\beta}(\log\log t_{n})^{\sigma+\bb\over 2\beta}\Big\}
<\infty\label{f5.1a}
\ee
for any $\th>1$ and 
\be
\sum_n \P\Big\{F(t_{n})
\ge \lambda_{2} t_{n}^{2\beta-\sigma\over 2\beta}
(\log\log t_{n})^{\sigma+\bb\over 2\beta}\Big\}
=\infty\label{f5.1a}
\ee
for $\th>1$ sufficiently large.

Conditioning on the stable process $\{X_t\}$, $F(t)$ is a centered
Gaussian process with variance
$$
\int_{\R^d}\xi(t,x)^2 dx.
$$
For $s<t$, the variance of $F(t)-F(s)$ is
$$
\int_{\R^d}\big[\xi(t,x)-\xi(s,x)\big]^2 dx.
$$
The relation
$$
\int_{\R^d}\big[\xi(t,x)-\xi(s,x)\big]^2 dx
\le \int_{\R^d}\xi(t,x)^2 dx-\int_{\R^d}\xi(s,x)^2 dx
$$
shows that, conditionally on on the stable process $\{X_t\}$, $F(t)$ is a $\P$-sub-additive Gaussian
process (\cite{KhoshnevisanLewis}). 
By Proposition 2.2, \cite{KhoshnevisanLewis},
$$
\P\Big\{\sup_{0\le s\le t}F(s)\ge a\Big\}\le 
2\P\Big\{F(t)\ge a\Big\},\hskip.2in a, t>0.
$$
Using (\ref{f5.1a}) and the Borel-Cantelli lemma we see that
\begin{equation}
\limsup_{n\to\infty}t_n^{-({2\beta-\sigma\over 2\beta})}(\log\log t_n)^{-\sigma+\bb\over 2\beta}
\sup_{t\leq t_{n}} F(t)\le\lambda_1,\hskip.2in a.s.\label{f5.2a}
\end{equation}
which leads to the desired upper bound in (\ref{5.6}).

It remains to establish the lower bound.

We will say that a countable set of  random variables $Y_{1}, Y_{2},\ldots$
is associated if for any $n$ and any coordinate-wise non-decreasing measurable functions
$f,g:R^{n}\mapsto R$ we have
\begin{equation}
\mbox{Cov }(f(Y_{1}, \cdots, Y_{n})g(Y_{1}, \cdots, Y_{n}))\geq 0,\label{9.1}
\end{equation}
and that the sequence (now the order counts) of  random variables $Y_{1}, Y_{2},\ldots$  is quasi-associated if for any $i<n$ and any coordinate-wise non-decreasing measurable functions
$f:R^{i}\mapsto R,\,\, g:R^{n-i}\mapsto R$ we have
\begin{equation}
\mbox{Cov }(f(Y_{1}, \cdots, Y_{i})g(Y_{i+1}, \cdots, Y_{n}))\geq 0.\label{9.2}
\end{equation}

Let $Y_{x},\,x\in Z^{d}$ be i.i.d. standard normals. By \cite[Theorem 2.1]{EPW} the set  
$Y_{x},\,x\in Z^{d}$ is associated. Let $h(x,y)$ be a non-negative function on $ Z^{d}\times Z^{d}$ such that
\begin{equation}
\wt h(x.x')=:\sum_{y\in Z^{d}}h(x,y)h(x',y)<\ff,\label{9.2a}
\end{equation}
and set
$V_{x}=\sum_{y\in Z^{d}}h(x,y)Y_{y}$. The collection $V_{x},\,x\in Z^{d}$ is  Gaussian process with covariance
\begin{equation}
E(V_{x}V_{x'})=f(x,x').\label{9.3}
\end{equation}
For each $m<\ff$, let $V_{m,x}=\sum_{y\in Z^{d},\,|y|\leq m}h(x,y)Y_{y}$. Since $V_{m,x}$ is a non-decreasing function of the $Y_{y},\,|y|\leq m$, it follows from
  \cite[$(P_{4})$]{EPW} that the set $V_{m,x},\,x\in Z^{d}$ is associated, and since by \cite[$(P_{5})$]{EPW} association is preserved under limits,  we also have that $V_{x},\,x\in Z^{d}$ is associated.

Let now $S=\{S_{0},S_{1},\ldots\}$ be a random walk in $Z^{d}$. Set $g_{h}(m)=\sum_{i=0}^{m}V_{S_{i}}$ and $g_{h}(k,l)=g_{h}(k)-g_{h}(l)$. It follows from the proof of 
\cite[Proposition 3.1]{KhoshnevisanLewis}  that for any $n$
and $0\leq s <t \leq u<v$,   
\begin{equation}
 g_{h}([ns],[nt])\hspace{.05 in}\mbox{ and }\hspace{.05 in} g_{h}([nu],[nv])\hspace{.05 in}\mbox{are quasi-associated.} \label{9.4}
\end{equation}

Using the stability of quasi-association under weak limits, we now show   that (\ref{9.4}) implies the following Lemma.
\begin{lemma}\label{lem-qa}
For any $0\leq s <t \leq u<v$, the pair 
\begin{equation}
F(t)-F(s), \hspace{.1 in}F(v)-F(u) \label{9.4a}
\end{equation}
is quasi-associated.
\end{lemma}

{\bf  Proof of Lemma \ref{lem-qa}: }

Note that we can write $g_{h}(m)=\sum_{y\in Z^{d}}\sum_{i=0}^{m}h(S_{i},y)Y_{y}$. Let
\begin{equation}
f_{\ep}(x)={e^{-\ep |x|} \over |x|^{ {\sigma +d\over 2}}+\ep}.\label{10.1}
\end{equation}
$f_{\ep}(x)$ is a positive, continuous integrable function of $x$, monotone decreasing in $|x|$.
We now define
\begin{equation}
h_{\ep,n,m}(x,y)={1 \over n\,m^{d/2}}f_{\ep}\({x \over n^{1/\bb}}-{y \over m}\).\label{10.2}
\end{equation}
It is clear that $h_{\ep,n,m}(x,y)$ satisfies (\ref{9.2a}). For notational conveinience we set 
\begin{equation}
g_{\ep,n,m}(r)=g_{h_{\ep,n,m}}(r).\label{10.2a}
\end{equation}
We take our random walk $S$ to be in the normal domain of attraction of $X$. We now show that
\be
\lim_{n\rar\ff}\(g_{\ep,n,m}([ns],[nt]),\,g_{\ep,n,m}([nu],[nv])\)\stackrel{w}{=}\(G_{\ep,m}(s,t),\,G_{\ep,m}(u,v)\), \label{10.3} 
\ee 
where
\begin{equation}
G_{\ep,m}(s,t)={1 \over m^{d/2}}\sum_{y\in Z^{d}}\(\int_{s}^{t}f_{\ep}(X_{r}-y/m)\,dr\)Y_{y}.\label{10.4}
\end{equation}
To see this we first note that
\begin{eqnarray}
&&E\Big (\exp \lc  i  \(ag_{\ep,n,m}([ns],[nt])+b\,g_{\ep,n,m}([nu],[nv])\Big)         \rc\)
\label{10.5}\\
&&=E\(e^{-H_{\ep,n,m}(a,b,s,t,u,v)/2}\),   \nonumber
\end{eqnarray}
where 
\begin{eqnarray}
H_{\ep,n,m}(a,b,s,t,u,v)&=&a^{2}\sum_{i,j=[ns]+1}^{[nt]}\wt h_{\ep, m} \({S_{i} \over n^{1/\bb}},{S_{j} \over n^{1/\bb}}\){1 \over n^{2}}
\label{10.6}\\
&&+2ab\sum_{i=[ns]+1}^{[nt]}\sum_{j=[nu]+1}^{[nv]}\wt h_{\ep, m} \({S_{i} \over n^{1/\bb}},{S_{j} \over n^{1/\bb}}\){1 \over n^{2}}\nonumber\\
&&+   b^{2}\sum_{i,j=[nu]+1}^{[nv]}\wt h_{\ep, m}  \({S_{i} \over n^{1/\bb}},{S_{j} \over n^{1/\bb}}\){1 \over n^{2}}\nonumber
\end{eqnarray}
and
\begin{equation}
\wt h_{\ep, m}(x,x')={1 \over  m^{d}}\sum_{y\in Z^{d}}f_{\ep}\(x-{y \over m}\) f_{\ep}\(x'-{y \over m}\).\label{10.6a}
\end{equation}
Similarly
\begin{eqnarray}
&&E\Big (\exp \lc  i  \(aG_{\ep,m}(s,t)+b\,G_{\ep,m}(u,v)\Big)         \rc\)
\label{10.5}\\
&&=E\(e^{-H'_{\ep,m}(a,b,s,t,u,v)/2}\),   \nonumber
\end{eqnarray}
where 
\begin{eqnarray}
&&H'_{\ep,m}(a,b,s,t,u,v)=a^{2}\int_{ [s,t]^{2}} \wt h_{\ep, m} (X_{r},X_{r'})\,dr\,dr'
\label{10.6}\\
&&+2ab\int_{ [s,t]\times [u,v]} \wt h_{\ep, m} (X_{r},X_{r'})\,dr\,dr'+   b^{2}\int_{ [u,v]^{2}} \wt h_{\ep, m} (X_{r},X_{r'})\,dr\,dr'.\nonumber
\end{eqnarray}
Thus to prove (\ref{10.3}) it suffices to show that
\be
\lim_{n\rar\ff}H_{\ep,n,m}(a,b,s,t,u,v)\stackrel{w}{=}H'_{\ep,m}(a,b,s,t,u,v), \label{10.3a} 
\ee 
and, since $\wt h_{\ep, m}(x,x')$ is continuous,  this follows directly from Skorohod's theorem, \cite{Sk}.

We next show that
\be
\lim_{m\rar\ff}\(G_{\ep,m}(s,t),\,G_{\ep,m}(u,v)\)\stackrel{w}{=}\(G_{\ep}(s,t),\,G_{\ep}(u,v)\), \label{10.3b} 
\ee 
where, using the notation of (\ref{f2}),
\begin{equation}
G_{\ep}(s,t)=\int_{R^{d}}\(\int_{s}^{t}f_{\ep}(X_{r}-x )\,dr\) \,W(dx) .\label{10.4b}
\ee 
But clearly
\begin{eqnarray}
&&E\Big (\exp \lc  i  \(aG_{\ep }(s,t)+b\,G_{\ep }(u,v)\Big)         \rc\)
\label{10.5b}\\
&&=E\(e^{-H'_{\ep }(a,b,s,t,u,v)/2}\),   \nonumber
\end{eqnarray}
where 
\begin{eqnarray}
&&H'_{\ep }(a,b,s,t,u,v)=a^{2}\int_{ [s,t]^{2}} \wt h_{\ep } (X_{r},X_{r'})\,dr\,dr'
\label{10.6b}\\
&&+2ab\int_{ [s,t]\times [u,v]} \wt h_{\ep } (X_{r},X_{r'})\,dr\,dr'+   b^{2}\int_{ [u,v]^{2}} \wt h_{\ep } (X_{r},X_{r'})\,dr\,dr'\nonumber
\end{eqnarray}
and now 
\begin{equation}
\wt h_{\ep}(x,x')= \int_{y\in R^{d}}f_{\ep}\(x-y\) f_{\ep}\(x'-y\)\,dy.\label{10.6b}
\end{equation}
Thus to obtain (\ref{10.3b}) it clearly suffices to show that
\be
\lim_{m\rar\ff}H'_{\ep,m }(a,b,s,t,u,v)\stackrel{^{L^{2}}}{=}H'_{\ep }(a,b,s,t,u,v).\label{10.3b} 
\ee 
To see this we note first that
\bea
\sup_{m,x,x'}\wt h_{\ep, m}(x,x') &\leq  &\sup_{m,x}{1 \over  m^{d}}\sum_{y\in Z^{d}}{e^{-\ep |x-y/m|} \over  \ep^{2}}\label{10.3b2}\\
&\leq  &\sup_{m,|x|\leq 1}{1 \over  m^{d}}\sum_{y\in Z^{d}}{e^{-\ep |x-y|/m} \over  \ep^{2}}\nn\\
&\leq  &\sup_{m}{1 \over  m^{d}}\sum_{y\in Z^{d}}{e^{-\ep |y|/m} \over  \ep^{2}}\leq C_{\ep},\nn
\eea
and therefore (\ref{10.3b}) follows easily from the dominated convergence theorem.
 
Finally, to prove our Lemma it now suffices to show that
\be
\lim_{m\rar\ff}\(G_{\ep}(s,t),\,G_{\ep}(u,v)\)\stackrel{w}{=}\(F(t)-F(s),\,F(v)-F(u)\). \label{10.3c} 
\ee 
As before, it  suffices  to show that
\be
\lim_{\ep\rar 0}H'_{\ep }(a,b,s,t,u,v)\stackrel{^{L^{2}}}{=}H'(a,b,s,t,u,v),\label{10.3c} 
\ee 
where 
\begin{eqnarray}
&&H' (a,b,s,t,u,v)=a^{2}\int_{ [s,t]^{2}} \wt h  (X_{r},X_{r'})\,dr\,dr'
\label{10.6c}\\
&&+2ab\int_{ [s,t]\times [u,v]} \wt h (X_{r},X_{r'})\,dr\,dr'+   b^{2}\int_{ [u,v]^{2}} \wt h  (X_{r},X_{r'})\,dr\,dr'\nonumber
\end{eqnarray}
and
\begin{equation}
\wt h (x,x')= \int_{\R^d}\vert x-y\vert^{-{\sigma +d\over 2}}\,\vert x'-y\vert^{-{\sigma +d\over 2}}\,dy
=c{1 \over |x-x'|^{\si}}.\label{10.6c}
\end{equation} 
Since  $\wt h_{\ep}(x,x')$ increases to $\wt h (x,x')$ as $\ep\rar 0$, (\ref{10.6c}) follows easily from Theorem \ref{eta-1} and the dominated convergence theorem.
\qed

To complete the proof of the lower bound, we recall from (\ref{eta-3a}) and (\ref{eta-3}) that
\bea
&&
E\((F(t)-F(s))(F(v)-F(u))\)\label{p1.1x}\\
&&
=c\E\int_{s}^{t}\!\!\int_{u}^{v}\vert X_a-X_b\vert^{-\sigma}\,da\,db\nonumber\\
&&
=c\E(\vert X_1\vert^{-\sigma})\,\,\int_{s}^{t}\!\!\int_{u}^{v}
{1\over  (a-b)^{\sigma/\beta}}\,da\,db.\nn
\eea
It then follows as in the proof of \cite[Theorem 5.2]{KhoshnevisanLewis} that for any $0<\la<1$,
if $0\leq s\leq \la t<t \leq u\leq \la v <v$, then 
\begin{equation}
\mbox{ Cov }\({F(t)-F(s) \over (t-s)^{1-\si/2\bb}},\,{ F(v)-F(u)\over (v-u)^{1-\si/2\bb}}\)
\leq c_{\la}\({t \over v}\)^{\si/2\bb}.\label{9.4b}
\end{equation}

The lower bound for our LIL then follows as in the proof of \cite[Theorem 1.1]{KhoshnevisanLewis}.
\qed

\section{Appendix: The limit as $M\to\ff$}\label{sec-Mlim} 
\medskip

\bt\label{theo-alt} Let $\rho_{\al,\ep}$ be defined in (\ref{32.1}) and
$\rho_{\alpha,\ep,M}$ be defined in (\ref{3.14}). We have
\begin{equation}
\limsup_{M\to\infty}M^{-d}\rho_{\alpha,\ep,M}\le \rho_{\al,\ep}.\label{ap.0}
\end{equation}
\et

\proof 
For any $x=(x_1,\cdots, x_d)\in\R^d$, we write
$[x]=([x_1],\cdots, [x_d])$ for the lattice part of $x$ (We also use the the notation
$[\cdots ]$ for parentheses without causing any confusion). For any $f\in {\cal
L}^2(\Z^d)$ with $\| f\|_2=1$,
\begin{eqnarray} &&
\sum_{\vert x\vert\le (2\pi)^{-1}Ma}  \wp_{\al, \ep}\({2\pi\over M}x\)
\bigg[\sum_{y\in\Z^d}\sqrt{Q\Big({2\pi\over M}(x+y)\Big)}
\sqrt{Q\Big({2\pi\over M}y\Big)}f(x+y)f(y)\bigg]^2\nn\\ &&
 =\int_{\{\vert \lambda\vert\le (2\pi)^{-1}Ma\}} \wp_{\al, \ep}\({2\pi\over
M}[\lambda]\)
\nn\\
&&\hspace{ .5in}
\bigg[\int_{\R^d}\sqrt{Q\Big({2\pi\over M}([\lambda]+[\gamma])\Big)}
\sqrt{Q\Big({2\pi\over M}[\gamma]\Big)}f([\lambda]+[\gamma])
f([\gamma])d\gamma\bigg]^pd\lambda\nonumber\\
 && =\Big({M\over
2\pi}\Big)^d\int_{\{\vert\lambda\vert\le a\}}
  \wp_{\al, \ep}\({2\pi\over M}[{M\over 2\pi}\lambda]\)\nonumber\\ &&
\hspace{1 in}\bigg[\Big({M\over 2\pi}\Big)^d
\int_{\R^d}\sqrt{Q_M
\Big(\gamma +{2\pi\over M}\Big[{M\over 2\pi}\lambda\Big]\Big)}
\sqrt{Q_M(\gamma)}  \nonumber\\ &&\hspace{1.5 in}\times f\Big(\Big[{M\over 2\pi}\lambda\Big]+
\Big[{M\over 2\pi}\gamma\Big]\Big) f\Big(\Big[{M\over
2\pi}\gamma\Big]\Big)d\gamma\bigg]^2d\lambda,\label{ap.14}
\end{eqnarray}
 where
\begin{equation} Q_M(\lambda)=Q\Big({2\pi\over M}\Big[{M\over
2\pi}\lambda\Big]\Big),
\hskip.2in\lambda\in\R^d.\label{64.15}
\end{equation}

Write
\begin{equation} g_0(\lambda) =\Big({M\over 2\pi}\Big)^{d/2}f\Big(\Big[{M\over
2\pi}\lambda\Big]\Big),
\hskip.2in\lambda\in\R^d.\label{ap.15}
\end{equation} We have
\begin{equation}
\int_{\R^d}g_0^2(\lambda)d\lambda=\Big({M\over 2\pi}\Big)^d
\int_{\R^d}f^2\Big(\Big[{M\over 2\pi}\lambda\Big]\Big)d\lambda
=\int_{\R^d}f^2([\lambda])d\lambda=\sum_{x\in\Z^d}f^2(x)=1.\label{}
\end{equation} We can also see that under this correspondence,
\begin{equation}
\Big({M\over 2\pi}\Big)^{d/2}f\Big(\Big[{M\over 2\pi}\lambda\Big]+
\Big[{M\over 2\pi}\gamma\Big]\Big)=g_0\Big(\gamma +{2\pi\over M}
\Big[{M\over
2\pi}\lambda\Big]\Big),\hskip.2in\lambda,\gamma\in\R^d.\label{ap.16}
\end{equation}
 
Therefore, we need only to show that for any fixed $a>0$
\begin{eqnarray} &&\qquad
\limsup_{M\to\infty}\sup_{\vert\vert g\vert\vert_2=1}
\int_{\{\vert \lambda\vert\le a\}} 
\wp_{\al, \ep}\({2\pi\over M}\Big[{M\over
2\pi}\lambda\Big]\)
\label{ap.17}\\ &&
\bigg[\int_{\R^d}\sqrt{Q_M
\Big(\gamma +{2\pi\over M}\Big[{M\over 2\pi}\lambda\Big]\Big)}
\sqrt{Q_M(\gamma)} g\Big(\gamma +{2\pi\over M}
\Big[{M\over 2\pi}\lambda\Big]\Big) g(\gamma)d\gamma\bigg]^2d\lambda
\nonumber\\ &&
 \le\sup_{\vert\vert g\vert\vert_2=1}\int_{\{\vert \lambda\vert\le a\}}
\hat{h}(\ep \la)\phi_{d-\si}( \la) 
\bigg[\int_{\R^d}\sqrt{Q(\lambda +\gamma)}
\sqrt{Q(\gamma)}g(\lambda +\gamma)g(\gamma)d\gamma\bigg]^2\,d\lambda.\nonumber
\end{eqnarray} 
To this end, note that by the inverse Fourier transformation  the function
\begin{equation} U_M(\lambda)=\int_{\R^d}\sqrt{Q_M (\gamma +\lambda)}
\sqrt{Q_M(\gamma)}g(\gamma +\lambda) g(\gamma)d\gamma\label{}
\end{equation} is the Fourier transform of the function
\begin{eqnarray} && V_M(x)={1\over
(2\pi)^d}\int_{\R^d}U_M(\lambda)e^{-i\lambda\cdot x}d\lambda\label{ap.18}\\
&&  ={1\over (2\pi)^d}\int_{\R^d}e^{-i\lambda\cdot
x}d\lambda\int_{\R^d}\sqrt{Q_M (\gamma +\lambda)}
\sqrt{Q_M(\gamma)}g(\gamma +\lambda) g(\gamma)d\gamma\nonumber\\ && 
={1\over (2\pi)^d}\int\!\!\int_{\R^d\times\R^d}e^{-i(\lambda-\gamma)\cdot x}
\sqrt{Q_M(\lambda)}g(\lambda)\sqrt{Q_M(\gamma)}g(\gamma)d\lambda
d\gamma\nonumber\\ &&  ={1\over
(2\pi)^d}\bigg\vert\int_{\R^d}e^{ix\cdot\gamma}\sqrt{Q_M(\gamma)}
g(\gamma)d\gamma\bigg\vert^2.\nonumber
\end{eqnarray} 
Therefore
\begin{eqnarray} &&
\int_{\R^d}\sqrt{Q_M
\Big(\gamma +{2\pi\over M}\Big[{M\over 2\pi}\lambda\Big]\Big)}
\sqrt{Q_M(\gamma)}g\Big(\gamma +{2\pi\over M}
\Big[{M\over 2\pi}\lambda\Big]\Big)g(\gamma)d\gamma \nn\\
&& =U_M\Big({2\pi\over M}\Big[{M\over 2\pi}\lambda\Big]\Big)\nn\\ &&
 ={1\over (2\pi)^d}\int_{\R^d}
\exp\Big\{ix\cdot {2\pi\over M}\Big[{M\over 2\pi}\lambda\Big]\Big\}
\bigg\vert\int_{\R^d}e^{ix\cdot\gamma}\sqrt{Q_M(\gamma)}
g(\gamma)d\gamma\bigg\vert^2dx\nonumber\\ &&
 \le {1\over (2\pi)^d}\int_{\R^d}\Big\vert 1-\exp\Big\{ix\cdot \Big(\lambda -
{2\pi\over M}\Big[{M\over 2\pi}\lambda\Big]\Big)\Big\}\Big\vert\nn\\ &&\hspace{1.5 in}\cdot
\bigg\vert\int_{\R^d}e^{ix\cdot\gamma}\sqrt{Q_M(\gamma)}
g(\gamma)d\gamma\bigg\vert^2dx\nonumber\\ &&\hspace{ 1in}
 +{1\over (2\pi)^d}\int_{\R^d} e^{ix\cdot \lambda}
\bigg\vert\int_{\R^d}e^{ix\cdot\gamma}\sqrt{Q_M(\gamma)}
g(\gamma)d\gamma\bigg\vert^2dx.\label{ap.19}
\end{eqnarray} 
By Parseval's identity and by the fact $Q_M\le 1$,
\begin{eqnarray} && {1\over (2\pi)^d}\int_{\R^d}\bigg\vert\int_{\R^d}
e^{ix\cdot\gamma}\sqrt{Q_M(\gamma)}
g(\gamma)d\lambda\bigg\vert^2dx\label{}\\ &&
=\int_{\R^d}Q_M(\gamma)g^2(\gamma)d\gamma
\le\int_{\R^d}g^2(\gamma)d\gamma =1.\nonumber
\end{eqnarray} Hence, the first term on the right hand side of (\ref{ap.19}) tends to
0 uniformly over
$\lambda\in\R^d$ and over all $g\in{\cal L}^2(\R^d)$  with $\vert\vert
g\vert\vert_2=1$ as $M\to\infty$.  The second term on the right hand side of
(\ref{ap.19}) is equal to
\begin{equation}
\int_{\R^d}e^{ix\cdot\lambda}V_M(x)dx=U_M(\lambda)=
\int_{\R^d}\sqrt{Q_M(\lambda +\gamma)}
\sqrt{Q_M(\gamma)}g(\lambda +\gamma)g(\gamma)d\gamma.\label{ap.21}
\end{equation}

Consequently, we will have (\ref{ap.17}) if we can prove
\begin{eqnarray} &&\quad
\limsup_{M\to\infty}\sup_{\vert\vert g\vert\vert_2=1}
\int_{\{\vert \lambda\vert\le a\}}  \wp_{\al, \ep}({2\pi\over M}[{M\over
2\pi}\lambda]) \label{ap.22}\\
&&\hspace{ 1in}\times
\bigg[\int_{\R^d}\sqrt{Q_M (\lambda +\gamma)}
\sqrt{Q_M(\gamma)}g(\lambda +\gamma)
g(\gamma)d\gamma\bigg]^2d\lambda\nn\\ &&\le\sup_{\vert\vert
g\vert\vert_2=1}\int_{\{\vert \lambda\vert\le a\}}\wp_{\al, \ep}(\lambda ) 
\bigg[\int_{\R^d}\sqrt{Q(\lambda +\gamma)}
\sqrt{Q(\gamma)}g(\lambda +\gamma)g(\gamma)d\gamma\bigg]^2d\lambda.
\nonumber
\end{eqnarray}

By uniform continuity of the function $Q$ we have that
$Q_M(\cdot)\to Q(\cdot)$ uniformly on $\R^d$. Thus, given $\epsilon >0$ we have
\begin{equation}
\sup_{\lambda,\gamma\in\R^d}\Big\vert\sqrt{Q_M (\lambda
+\gamma)}\sqrt{Q_M(\gamma)}-\sqrt{Q (\lambda +\gamma)}
\sqrt{Q(\gamma)}\Big\vert <\epsilon\label{ap.23}
\end{equation} for sufficiently large $M$. Therefore,
\begin{eqnarray} &&
\bigg\{\int_{\{\vert \lambda\vert\le a\}}d\lambda
\bigg[\int_{\R^d}\sqrt{Q_M (\lambda +\gamma)}
\sqrt{Q_M(\gamma)}g(\lambda +\gamma)
g(\gamma)d\gamma\bigg]^2\bigg\}^{1/2}\label{ap.24}\\ &&
\le \epsilon\bigg\{\int_{\{\vert \lambda\vert\le a\}}d\lambda
\bigg[\int_{\R^d}g(\lambda +\gamma)
g(\gamma)d\gamma\bigg]^2\bigg\}^{1/p} \nonumber
\\ && +\bigg\{\int_{\{\vert \lambda\vert\le a\}}d\lambda
\bigg[\int_{\R^d}\sqrt{Q (\lambda +\gamma)}
\sqrt{Q(\gamma)}g(\lambda +\gamma)
g(\gamma)d\gamma\bigg]^2\bigg\}^{1/2}.\nonumber
\end{eqnarray} 
Also, since $\|g\|_{2}=1$,
\begin{equation}
\int_{\{\vert \lambda\vert\le a\}}d\lambda
\bigg[\int_{\R^d}g(\lambda +\gamma) g(\gamma)d\gamma\bigg]^2\le
C_da^d,\label{ap.25}
\end{equation} where $C_d$ is the volume of a $d$-dimensional unit ball. (\ref{ap.22}) then follows using the uniform continuity of $\wp_{\al,\ep}(\lambda )$.
\qed

{\bf  Acknowledgement: } We would like to thank Frank den Hollander for useful information about polymer models.

\end{document}